\documentclass[final]{amsart}

\newcommand{\SubjClass}[2]{\subjclass[2000]{Primary #1; Secondary #2}}
\newcommand{\Title}[2][]{\title[#1]{#2}}

\makeindex
\usepackage[notref,notcite]{showkeys}

\usepackage{amssymb}
\providecommand{\qedhere}{\relax}

\usepackage[all]{onlyamsmath}
\usepackage{mathrsfs}
\usepackage{xspace}
\usepackage{ifthen}
\newcommand{\Ifempty}[3]{\ifthenelse{\equal{#1}{}}{#2}{#3}}
\newcommand{\Ifnempty}[2]{\Ifempty{#1}{}{#2}}

\usepackage{verbatim}
\usepackage[matrix,arrow]{xy}

\usepackage{enumerate}
\usepackage[bookmarks,breaklinks,colorlinks,unicode]{hyperref}

\newcommand{\newthm}[2]{\newtheorem{#1}[lemma]{#2}}
\providecommand{\newnumbered}[1]{\newtheorem{#1}}
\newcommand{\newnumb}[2]{\newnumbered{#1}[lemma]{#2}}
\theoremstyle{plain}
\newtheorem{lemma}{Lemma}
\newnumb{Ass}{Assumption}
\newnumb{convention}{Convention}
\newnumb{criterion}{Criterion}
\newnumb{fact}{Fact}
\newthm{cor}{Corollary}
\newthm{claim}{Claim}
\newthm{thm}{Theorem}
\newthm{prop}{Proposition}
\newnumb{axiom}{Axiom}

\theoremstyle{definition}
\newnumb{definition}{Definition}
\newnumb{exampl}{Example}

\newnumb{remrk}{Remark}
\newnumb{notation}{Notation}

\newenvironment{example}{\begin{exampl}}{\hspace{\stretch{1}}\rule{1ex}{1ex}\end{exampl}}
\newenvironment{remark}{\begin{remrk}}{\end{remrk}}

\newcommand{\Def}[1]{\emph{#1}\xspace}

\newcommand{\DefAlias}[2]{\expandafter\xdef\csname #1\endcsname{#2}}
\newcommand{\CiteAlias}[2]{\DefAlias{CITE#1}{#2}}
\newcommand{\Cite}[2][]{\Ifempty{#1}{\cite{\csname CITE#2\endcsname}}{\cite[#1]{
\csname CITE#2\endcsname}}}

\newcommand{\Zilber}{Zil$'$ber}

\CiteAlias{ACFA}{MR1652269}
\CiteAlias{LAG}{MR1102012}
\CiteAlias{DG}{MR1972449}
\CiteAlias{groupoid}{groupoid}
\CiteAlias{bounds}{MR739626}
\CiteAlias{stability}{MR1416106}
\CiteAlias{SGA}{MR0354652}
\CiteAlias{singer}{MR1480919}
\CiteAlias{milne}{MR559531}
\CiteAlias{topoi}{MR1300636}
\CiteAlias{sacks}{MR0398817}
\CiteAlias{EI}{MR1083551}
\CiteAlias{HrWeil}{MR1827833}
\CiteAlias{SingerMod}{MR2441376}
\CiteAlias{PillayZoe}{MR1650667}
\CiteAlias{Bouscaren}{MR984628}
\CiteAlias{eisenbud}{MR1322960}
\CiteAlias{ZilBind}{MR606796}
\CiteAlias{Poizat}{MR727805}
\CiteAlias{PillayDCFI}{MR1649893}
\CiteAlias{PillayDGPS}{MR2094550}
\CiteAlias{HartShami}{MR2140036}
\CiteAlias{Pezz}{MR1919592}
\CiteAlias{SingerDL}{MR1960772}
\CiteAlias{Kolchin}{MR0568864}
\CiteAlias{Deligne}{MR1106898}
\CiteAlias{Etingof}{MR1336694}
\CiteAlias{Franke}{MR0155819}
\CiteAlias{ogroups}{MR1081816}
\CiteAlias{Pillay}{MR1429864}
\CiteAlias{prodef}{MR2335085}

\newcommand{\Email}{\url{mailto:mkamensky@math.uwaterloo.ca}}

\newcommand{\x}{\times}
\newcommand{\MapsTo}{\rightarrow}
\newcommand{\df}{\stackrel{\mathrm{def}}{=}}
\newcommand{\w}{\mbox{${\omega}$}}
\newcommand{\ti}[1]{\tilde{#1}}
\newcommand{\bigand}[1]{{\bigwedge}_{#1}}
\newcommand{\tensor}{\otimes}

\newcommand{\Emph}[1]{\emph{#1}\index{#1}}

\newcommand{\QQ}{\mathbb{Q}}

\newcommand{\Pro}[2][]{\operatorname{Pro}_{#1}({#2})}
\newcommand{\Ind}[2][]{\operatorname{Ind}_{#1}({#2})}
\newcommand{\Eq}[1]{{\{#1\}}^{eq}}
\newcommand{\EQ}[1]{{#1}^{eq}}
\newcommand{\Q}[1][]{\mathbf{Q}\Ifnempty{#1}{(#1)}}
\newcommand{\Cc}[1][]{\mathbf{C}\Ifnempty{#1}{(#1)}}
\newcommand{\X}[1][]{\mathbf{X}\Ifnempty{#1}{(#1)}}
\newcommand{\DCF}{$DCF$}
\newcommand{\ACFA}{$ACFA$}
\newcommand{\ACF}{$ACF$}
\newcommand{\sg}{\sigma}
\newcommand{\ogroup}{$\w$-group\xspace}
\newcommand{\ogroups}{$\w$-groups\xspace}
\newcommand{\Mat}[4]{\bigl(\begin{smallmatrix}{#1}&{#2}\\{#3}&{#4}\end{smallmatrix}\bigr)}
\newcommand{\Dcl}[1][]{dcl\Ifnempty{#1}{(#1)}}
\newcommand{\Tp}[2]{\mbox{$tp(#1/_{#2})$}}
\newcommand{\Gl}[2][]{GL_{#2\Ifnempty{#1}{,}#1}}
\newcommand{\ra}{\rightarrow}
\newcommand{\onto}{\twoheadrightarrow}
\newcommand{\autQ}[2]{\inv{f_{#2}}\circ{f_{#1}}}
\newcommand{\autC}[2]{{f_{#2}}\circ\inv{f_{#1}}}
\newcommand{\inv}[1]{{#1}^{-1}}
\newcommand{\Conj}{\sim}
\newcommand{\Inf}{\infty}
\newcommand{\QF}{\Delta}
\newcommand{\Aut}[2][]{Aut_{\QF}(#2\Ifnempty{#1}{/#1})}
\newcommand{\Z}{\mathbb{Z}}
\newcommand{\Se}[2][]{{#2}^{SE}_{#1}}
\newcommand{\LL}{\mathcal{L}}
\newcommand{\ka}{\kappa}
\newcommand{\Ob}{\mathfrak{Ob}}
\newcommand{\Mor}{\mathfrak{Mor}}
\providecommand{\G}{}
\renewcommand{\G}{\mathcal{G}}
\newcommand{\isom}{\stackrel{\sim}{\rightarrow}}
\newcommand{\tp}[2][C]{\mbox{$tp_{\scriptscriptstyle\QF}(#2/_{#1})$}}
\newcommand{\Acl}[1][]{acl\Ifnempty{#1}{(#1)}}
\newcommand{\Df}[2][]{\mathcal{D}_{#1}(#2)}
\newcommand{\Stn}[2][]{\mathcal{S}_{#1}(#2)}
\newcommand{\Par}[2]{#1|_{#2}}
\newcommand{\Res}[2]{#1\upharpoonright_{#2}}
\newcommand{\Cb}[2][]{\mathscr{C}b_{#1}(#2)}
\newcommand{\kk}{\Bbbk}
\newcommand{\Singer}[1]{#1 in \Cite{singer}}
\newcommand{\DHom}[3][\kk]{Hom^{\sg}_{#1}(#2,#3)}
\newcommand{\Det}[1]{\det(#1)}
\newcommand{\A}{\mbox{$\mathbb A$}}
\newcommand{\subgrpeq}{\leq}
\newcommand{\Wbar}[1]{\overline{#1}}
\newcommand{\what}[1]{\widehat{#1}}

\Title[definable automorphism groups]{Definable groups of partial automorphisms}
\author{Moshe Kamensky}
\SubjClass{03C40}{03C60}

\renewcommand{\G}{\mathcal{G}}
\newenvironment{Thanks}{\subsection*{Acknowledgement}}{}
\newenvironment{Proof}[1][]{\begin{proof}[Proof #1]}{\end{proof}}
\newcommand{\MakePreTitle}{}
\newcommand{\MakePostTitle}{\maketitle}

\address{Department of Mathematics\\
         The Hebrew University\\
         Jerusalem, Israel}
\curraddr{Department of pure Math\\
         University of Waterloo\\
         Waterloo, On, N2L 3G1, Canada}

\email{\Email}
\urladdr{\url{http://mkamensky.notlong.com}}

\begin{document}


\MakePreTitle
\begin{abstract}
  The motivation for this paper is to extend the known model theoretic 
  treatment of differential Galois theory to the case of linear 
  \emph{difference} equations (where the derivative is replaced by an 
  automorphism.) The model theoretic difficulties in this case arise from the 
  fact that the corresponding theory ACFA does not eliminate quantifiers. We 
  therefore study groups of restricted automorphisms, preserving only part of 
  the structure. We give conditions for such a group to be (infinitely) 
  definable, and when these conditions are satisfied we describe the 
  definition of the group and the action explicitly.

  We then examine the special case when the theory in question is obtained by 
  enriching a stable theory with a generic automorphism. Finally, we 
  interpret the results in the case of ACFA, and explain the connection of 
  our construction with the algebraic theory of Picard-Vessiot extensions.

  The only model theoretic background assumed is the notion of a definable 
  set.
\end{abstract}

\MakePostTitle

\section{Introduction}\label{sec:intro}
A \emph{linear differential equation} is an equation of the form \(Dx=Ax\), 
where \(D\) is a (formal) derivation, \(A\) is a matrix over some base 
differential field, and \(x\) is a tuple of variables. To any such equation, 
it is possible to associate a certain extension of differential fields, the 
\Emph{Picard-Vessiot extension}, that contains a system of solutions to this 
equation. The \Emph{Galois group} of the equation is defined to be the 
automorphism group of this field. In~\Cite{DG}, it is shown that when the 
base field is \(\QQ(t)\), with \(Dt=1\), this Galois group is always 
computable.

The fundamental observation for the results of that paper, is that there is a 
model theoretic interpretation of this Galois group. More precisely, there is 
a general definition of the notion of ``the group of automorphisms of a 
definable set \(\Q\) over another definable set \(\Cc\)''\index{automorphisms 
group|(}. When \(\Q\) is internal to \(\Cc\) (i.e., has a definable family of 
bijections into \(\Cc\); see section~\ref{ssc:Internality} for the definition), 
this group turns out to be the group of points of a type-definable group.  
This fact is also explained in~\Cite{DG}, in appendix B.

To apply the general construction to linear differential equations, one 
considers the theory of differentially closed fields (\DCF{}), the model 
completion of the theory of differential fields (with constant symbols for 
the base field.) The equation \(Dx=Ax\) is then a interpreted as the definable 
set \(\Q\), while the set of constants \(Dx=0\) plays the role of \(\Cc\). The 
internality condition corresponds to the fact the \(\Q\) is a finite 
dimensional vector space over \(\Cc\), and thus has a definable family of 
bijections with some power of \(\Cc\). To identify the model theoretic group 
with the (algebraic) Galois group, one may embed the Picard-Vessiot extension 
into a model of \DCF{}, and then use the fact that \DCF{} eliminates 
quantifiers.

The purpose of the present work is to describe in more detail the 
construction of the model theoretic group of automorphisms, and to generalise 
it in a way that will be suitable for dealing with \emph{difference} 
equations. A linear difference equation\index{linear difference equation} is 
an equation of the form \(\sg(x)=Ax\), where \(\sg\) is a formal automorphism, 
and \(A\), as before, is a matrix over some base difference field (i.e., a 
field with a prescribed automorphism.) The algebraic theory of this case is 
described in~\Cite{singer}. It is analogous to the case of differential 
equations, but differs in some points.  In particular, the Galois group is 
only constructed in the case that the subfield of constant elements of the 
base field is algebraically closed.

From the model theoretic point of view, there are two essential differences 
between this case and the case of differential equations: First, the theory 
of algebraically closed fields with an automorphism (\ACFA{}), which is the 
analogue of \DCF{} in this case, does not eliminate quantifiers. In 
particular, there are definable subsets of the set defined by the equation 
that are not algebraic (quantifier free definable), and the original 
construction would produce the group of automorphisms preserving all such 
definable sets, not only the quantifier free ones (see also 
example~\ref{xpl:noextend}, proposition~\ref{prp:profinite} and 
example~\ref{xpl:profinite}).  The algebraic Galois group, in contrast, 
preserves the algebraic structure only.  Thus we need to construct a group of 
automorphisms that preserves only some definable sets.

The other distinction of this case is that, though the Galois group is 
constructed as the automorphism group of a certain Picard-Vessiot extension,
this extension may have zero-divisors. Therefore, even when the quantifier 
free group is constructed, it is not clear that it coincides with the 
algebraic Galois group.

There are also intrinsic questions associated with the original 
constructions. Specifically, if the automorphism group is definable, we may 
consider its points in an arbitrary model, not just a saturated one. Our new 
description interprets every such group as a group of automorphisms. 
Additionally, the formulas defining the group are produced explicitly.

We now briefly summarise the contents that follow.  In 
section~\ref{sec:general} we give the definitions of automorphism groups and 
of internality, and describe conditions for the automorphism group \(G\) to 
be (type-) definable. In the case when these conditions are satisfied, we 
show this by presenting the definition of the group explicitly. We also 
compare this group with the original construction in~\Cite{DG}. We consider 
the dependence of the automorphism group on the internality datum, and show 
that it is essentially independent. Finally, we describe a definable family 
of groups, that act on the \(G\)-torsor used to construct \(G\). The 
algebraic Galois group, considered in~\Cite{singer} and in 
section~\ref{sec:ACFA}, is eventually identified with a member of the Zariski 
closure of this family.

In section~\ref{sec:stability} we approach closer to the example of \ACFA{}, 
and consider theories obtained from a stable theory by adding a generic 
automorphism. The goal is to obtain a more precise description of the 
definition of the group, as given by equations over the base structure. In 
the case of \ACFA{} (where the stable theory is \ACF{}), this means polynomial 
equations. It turns out that the existence of such a description follows 
essentially from the stability alone.

Finally, in section~\ref{sec:ACFA} we consider the case of \ACFA{} itself. We 
describe the interpretation of the structure we obtained in previous sections 
for this case, as well as some more specific features that follow mainly from 
the Noetherian property in this case. We then explain the connection with the 
algebraic Galois group of~\Cite{singer}, and also consider some examples.
\index{automorphisms group|)}

For completeness, an appendix is included where the basic model-theoretic 
results used in the paper are briefly explained. This should hopefully make 
the paper accessible to readers with no model theoretic background.

\begin{Thanks}
This work is part of my PhD research, performed in the Hebrew university 
under the supervision of Ehud Hrushovski. I would like to thank him for his 
guidance, and in particular for suggesting this subject, which, among other 
things helped me to begin understanding stability theory.

I would also like to thank the referee for the careful reading and useful 
comments.
\end{Thanks}

\subsection{History}
Definable automorphism groups first appeared in the context of strongly 
minimal theories, in works of \Zilber{} (cf.~\Cite{ZilBind}), where they are 
called \Emph{binding groups}. Poizat observed, in~\Cite{Poizat}, that this 
construction can be used to give a model theoretic interpretation of 
differential Galois theory. This approach was taken further by Pillay, 
in~\Cite{PillayDCFI} (and~\Cite{PillayDGPS}), who extended the algebraic work 
of Kolchin.

In the abstract setting, the results on the binding groups were extended by 
Hrushovski to the stable case in~\Cite{ogroups}. They were further extended 
to simple theories in various places, including~\Cite{HartShami} 
and~\Cite{Pezz}. The ultimate result, on which this work is based, appears in 
appendix B of~\Cite{DG}, where the automorphism group is proved to exist with 
the sole assumption of internality. The main addition of the current work in 
this respect is the elementary and explicit description of this group, from 
which the extension to partial automorphisms is clear (as well as the fact 
that the binding group is an intersection of definable groups, even in the 
unstable case.) The interpretation of this theory in terms of definable 
groupoids\index{groupoid!definable} appears in~\Cite{groupoid}.

Theories of the form \(T_\sg\)\index{aa@\(T_\sg\)} appearing in 
section~\ref{sec:stability} were considered by Pillay and Chatzidakis 
in~\Cite{PillayZoe} (claim~\ref{clm:acl}, for example, is proved there.) 
Elimination of imaginaries for such theories is characterised (in terms of 
\(T\)) in~\Cite{groupoid}.

The motivating example for a big part of this theory is differential Galois 
theory. As mentioned above, the model theoretic connection was observed by 
Poizat and further developed by Pillay (and by Hrushovski in~\Cite{DG}.) The 
algebraic theory was known in special cases to 19th century mathematicians, 
and was systematically developed by Kolchin (\Cite{Kolchin}), and later by 
many authors, of which Singer and van-der-Put (\Cite{SingerDL}) and Deligne 
(\Cite{Deligne}) seems the most relevant to the model theoretic approach.

In the case of \emph{difference} equations, the theory seems to be much less 
developed. The algebraic theory, as we treat it, is developed by Singer and 
van-der-Put in~\Cite{singer}. A slightly different approach is taken 
in~\Cite{Franke} and in~\Cite{Etingof}. For the model theoretic context, the 
theory we use was developed by Hrushovski and Chatzidakis in~\Cite{ACFA}. The 
model theory of the Galois group was recently considered in~\Cite{SingerMod}.

\subsection{Notation and conventions}\label{ssc:notation}
We consider an arbitrary (not necessarily complete) theory \(T\). By a 
\Def{definable set} we mean a formula in the language of \(T\), up to 
equivalence with respect to \(T\) (two formulas are equivalent if the define 
the same subset in any model of \(T\)). In particular, definable sets are 
always over \(0\) unless explicitly mentioned otherwise.

If \(A\) is any subset of a model \(M\), \(\Dcl(A)\)\index{aa@\(\Dcl(A)\)} is 
the set of all elements in \(M\) of the form \(f(\Wbar{a})\), where 
\(\Wbar{a}\) is a tuple of elements from \(A\), and \(f\) is a 
\(T\)-definable function on a definable set containing \(\Wbar{a}\).Note that 
this set does not depend on the ambient model \(M\), but only on its theory 
(in the language with constant symbols for \(A\)): \(\Dcl(A)\) computed in 
one model of this theory is in canonical bijection over \(A\) with 
\(\Dcl(A)\) computed in any other. Given a definable (or, more generally, 
ind-definable) set \(X\), we  denote by \(X(A)\) the set of points of \(X\) 
that belong to \(\Dcl(A)\).

We will use the notion of a \Emph{definable family of (definable) sets}. A 
family is simply a definable set \(\phi(x,y)\subseteq{}X\x{}Y\) where one of 
the variables, say \(x\), is considered the parameter variable of the family.  
Thus, in this case we have a family of subsets of \(Y\), varying with \(x\).  
In many cases we will require the members of the family (the fibres) to be 
disjoint. A general family \(\phi(x,y)\) can be converted to a disjoint one 
with the same parameter variable and isomorphic fibres via the formula 
\(\forall{}x,y,z(\psi(x,z,y)\iff(x=z\land\phi(x,y)))\) (which may require a 
quantifier).

In appendix~\ref{sec:app}, we recall some other basic notions of model 
theory, namely elimination of imaginaries and stable embeddedness. One 
unusual aspect of our setting with respect to these notions is that we do not 
assume \(T\) to be complete. This requires some extra care with the 
definitions, but does not present real difficulties. This is also explained 
in the appendix.

We shall mildly use the notions of pro-definable and ind-definable sets. Some 
details on this can be found in~\Cite{prodef}. For simplicity, one can always 
think about a pro-definable set as a partial type, and of an ind-definable 
set as a bounded increasing union of definable sets. For groups, we use the 
following terminology:
\begin{definition}
  An \Def{\ogroup}\index{group@\ogroup} in a theory \(T\) is the intersection 
  of a (language sized) chain of definable groups
\end{definition}
We note that in stable theories, any pro-definable group has this form 
(cf~\Cite[Lemma~6.18]{Pillay}), but this is false in general.

\section{Internality and definable Galois groups}\label{sec:general}
In this section, we deal with the basic notion of a set \(\Q\) being 
\emph{internal} to \(\Cc\). This roughly means that \(\Q\) has a definable 
family of bijections with \(\Cc\). We shall see that having this situation is 
almost equivalent to automorphism groups of \(\Q\) over \(\Cc\) being 
\ogroups\index{group@\ogroup}.

A basic example is as follows: Let \(\Q\) be a vector space of dimension \(n\) 
over a field, and let \(\Cc\) the \(n\)-th Cartesian power of the field. The 
family of vector space bijections between \(\Q\) and \(\Cc\) can be identified 
with the set \(\X\) of (ordered) vector space bases of \(\Q\). In this 
identification, an element \(x\) of \(\X\) maps a vector in \(\Q\) to the 
coefficients of its presentation in that basis. This set \(\X\), as well as the 
family of maps from \(\Q\) to \(\Cc\) are definable from the vector space 
structure. Any element of the group of linear automorphisms \(\Gl{}(\Q)\) of 
\(\Q\) can be obtained as a composition of one bijection of this type, with the 
inverse of another. Automorphism groups preserving any additional structure 
will be definable sub-groups (or \ogroups) of this group.

In general, given arbitrary sets \(\Q\) and \(\Cc\), and a set \(\X\) of 
injective maps from \(\Q\) to \(\Cc\), consider the set \(F\) of bijections of 
\(\Q\) to itself of the form \(\inv{f}\circ{}g\), with \(f,g\in{}\X\) (for 
which the composition is defined.) Composing such a bijection \(t\in{}F\) 
with an element \(h\in\X\) gives a new function \(h\circ{}t\) from \(\Q\) to 
\(\Cc\). If, for some \(t\), this new function is again in \(\X\) for any 
\(h\in\X\), we say that \(t\) preserves \(\X\).  The set of elements of \(F\) 
preserving \(\X\) in this manner forms a group \(G\).  This group acts on 
\(\Q\) by evaluation and on \(\X\) by composition, and turns out to be the 
group of automorphisms of the evaluation map from \(\Q\x{}\X\) to \(\Cc\).  
When all of the above structure arises as the sets of points of definable 
sets, \(\Q\) is said to be \emph{internal} to \(\Cc\), and the above 
construction is carried out definably in subsection~\ref{ssc:Internality}.  
As in the case of vector spaces, any additional (definable) structure is the 
preserved by a definable subgroup of this group.

Internality is defined in definition~\ref{dfn:Internality}. Automorphisms and 
groups of them are defined in section~\ref{ssc:partial}. The final result is 
stated as theorem~\ref{thm:abstract}. Throughout the section, unless 
mentioned otherwise, we work with arbitrary models of an arbitrary theory 
\(T\), which is not necessarily complete, but which eliminates imaginaries 
(possibly by passing to \(T^{eq}\).

\subsection{Internality}\label{ssc:Internality}
\index{internality|(}
We are interested in the following situation:

\begin{definition}\label{dfn:Internality}
  Let \(\Q\), \(\Cc\) be definable sets.
  \begin{enumerate}
    \item
      \(\Q\) is said to be 
      \emph{\((\X,f)\)-internal}\index{internal!\((\X,f)\)-} (or simply 
      \Def{internal}) to \(\Cc\) if \(\X\) is a definable set and 
      \(f:\Q\x\X\ra\Eq{\Cc}\) is a definable map, such that, for any 
      \(x\in\X\), the map \(f_x:\Q\ra\Eq{\Cc}\) defined as \(f_x(q)=f(q,x)\) 
      is injective, and for \(x\ne{}y\), \(f_x\) and \(f_y\) are distinct.
      
    \item
      If \(\Q\) is \((\X,f)\)-internal to \(\Cc\), and \(M\) is a model, we 
      denote by 
      \[Aut_{(\X,f)}(\Q/\Cc)(M)\]\index{aa@\(Aut_{(\X,f)}(\Q/\Cc)(M)\)} the 
      group of pairs \((\tau_{\Q},\tau_{\X})\), where 
      \(\tau_{\Q}:\Q(M)\ra\Q(M)\) and \(\tau_{\X}:\X(M)\ra\X(M)\) are 
      bijections, and
      \begin{equation*}
        f(\tau_{\Q}(q),\tau_{\X}(x))=f(q,x)
      \end{equation*}
      for any \((q,x)\in\Q\x\X\). This group will be called the 
      \Def{internality group}.

    \item
      Assume that for any model \(M\), we are given a group \(G(M)\) acting 
      faithfully on \(\Q(M)\) (\(G\) is not assumed to be functorial in 
      \(M\).) \emph{\(\Q\) is internal to \(\Cc\) relatively to 
      \(G\)}\index{internal!relatively to \(G\)} if \(\Q\) is 
      \((\X,f)\)-internal to \(\Cc\) for some \(\X,f\) such that, for any 
      model \(M\), any \(g\in{}G(M)\) extends to an element of 
      \(Aut_{(\X,f)}(\Q/\Cc)(M)\).

      In this case, the pair \((\X,f)\) is called an \emph{internality 
      datum}\index{internality!datum} for \(G\).

  \end{enumerate}
\end{definition}
\begin{remark}\label{rmk:Internality}
  \mbox{}
\begin{enumerate}
\item
  An element of \(G(M)\) is thought of as an automorphisms of \(\Q(M)\).  The 
  condition of internality relatively to \(G\) says that any such automorphism 
  can be extended to an automorphism of the internality structure. Our goal 
  is to find conditions that \(G\) is an \ogroup (in particular, the data is 
  only interesting when the isomorphism class of \(G(M)\) depends only on the 
  isomorphism class of \(M\).) Note that in this part of the definition, 
  \((\X,f)\) is not part of the data, and \(G\) does not depend on them.

\item \label{srk:distinctn}
  Given a family \((\X,f)\) as in the definition of internality, except that 
  the \(f_x\) are not distinct, we may obtain an internality datum from it 
  via elimination of imaginaries, with the same functions (from \(\Q\) to 
  \(\Cc\)) appearing as fibres.

\item
  If \(\Q\) is \(\Cc\) internal relatively to \(G\), then the same is true 
  for any \(G_1\subgrpeq{}G\) (i.e., \(G_1(M)\) is a subgroup of \(G(M)\), 
  for any \(M\).) In other words, if \(\Q\) is \((\X,f)\)-internal to 
  \(\Cc\), then \(\Q\) is \(\Cc\) internal relatively to any 
  \(G\subgrpeq{}Aut_{(\X,f)}(\Q/\Cc)\).

\item
  Internality can be defined in the same way when \(\Cc\) is 
  ind-definable\index{ind-definable}, rather than just definable. However, it 
  amounts to saying the same thing for some definable subset of \(\Cc\) in this 
  case, and therefore for simplicity we don't do it.

\item
  As can be seen from the definition (and below), it is harmless to replace 
  \(\Cc\) by \(\Cc'\) if \(\Eq{\Cc}=\Eq{\Cc'}\), and, in fact, by any 
  \(\Cc'\) such that \(\Eq{\Cc'}\subseteq{}\Eq{\Cc}\) and contains sorts for 
  all definable sets in question. We shall assume, below, that the values of 
  \(f\) lie, in fact, in \(\Cc\).
\end{enumerate}
\end{remark}

Note that the definition does not require that individual members of the  
automorphisms group are definable (over any parameters.) However, this is in 
fact the case, as follows from proposition~\ref{prp:autispair} below.

In the following parts we will consider some auxiliary definable sets that 
appear as a result of internality. The reader may find it convenient to refer 
to the summary of notation in~\ref{ssc:summary}.

\subsubsection{The set of images}\label{sss:images}
Let \(\Q\) be \((\X,f)\) internal to \(\Cc\). Different elements of \(\X\) 
may map \(\Q\) to distinct subsets of \(\Cc\). The family of all these 
subsets is parametrised by a definable (in \(\EQ{T}\)) set \(D\). We apply to 
this family the process described in~\ref{ssc:notation}, namely, we replace 
\(\Cc\) by \(\Cc\x{}D\) (denoting this new set by \(\Cc\) again), so that the 
family of image sets is given as fibres of a projection to \(D\).  The 
internality datum is modified accordingly.  We have a natural map 
\(\pi:\X\onto{}D\), sending each element of \(x\in\X\) to the image of \(\Q\) 
under \(x\).  The automorphism group of the internality structure preserves 
this map: For any \(\tau\in{}Aut_{(\X,f)}(\Q/\Cc)(M)\), 
\(\pi(\tau(x))=\pi(x)\) (indeed, \(f(\tau(q),\tau(x))=f(q,x)\), so \(f_x\) 
and \(f_{\tau(x)}\) have the same image.)

Thus, the internality datum can be described as follows: We are given a 
family of maps \(f:\Q\x\X\ra\Cc\), another set \(D\), and maps from \(\X\) 
and from \(\Cc\) to \(D\) (whose fibres over a point \(d\in{}D\) will be 
denoted \(\X_d\) and \(\Cc_d\)), such that the combined map 
\((f,p_2):\Q\x\X\ra\Cc\x_D\X\) (where \(p_2\) is the second projection) is a 
bijection. From now on the internality datum will be assumed to be in this 
form. The map \(\Cc\x_D\X\ra\Q\) obtained from the inverse of the above map 
will be denoted by \(g\). As mentioned above, if this data witnesses 
internality relatively to \(G\), the action is naturally induced on 
\(\Cc\x_D\X\), and we have \(g(c,h(x))=h(g(c,x))\) for any \(h\in{}G\).  If 
all elements of \(\X\) map \(\Q\) to the same set, \(D\) is one point, and we 
have a family of bijections between \(\Q\) and \(\Cc\). For our purposes, 
there is no substantial difference between this case and a general \(D\). At 
the other extreme, if all the image sets are different, then in the new 
description \(\Q\x\X\) has a definable bijection with \(\Cc\). In particular, 
\(G\) is trivial, and if \(\Cc\) is stably embedded, then \(\Q\) is a subset 
of \(\Eq{\Cc}\).

\begin{remark}\label{rmk:oldint}
  In appendix B of~\Cite{DG}, internality is defined in terms of \(g\) rather 
  than \(f\), and \(g\) is not required to be injective. In other words, 
  according to that definition, \(\Q\) is internal to \(\Cc\) if there is a 
  definable map \(g:\Cc\x_D\X\ra\Q\) with each \(g_x\) surjective. The kernel 
  of \(g\) is then an equivalence relation on \(\Cc\), and dividing by it, we 
  get a new set in \(\Eq{\Cc}\), where the induced map is bijective, on the 
  fibres.  We thus get back to our situation. Since any group of automorphism 
  considered (either here or in~\Cite{DG}) acts trivially on \(\Eq{\Cc}\), 
  this process does not affect the automorphism group.

  Equivalently, instead of dividing by the kernel of \(g\), we may work with 
  the original datum, but modify \(D\) to be the family of domains of \(g\), 
  together with the equivalence relation given by the kernel of \(g\).

  See also~\ref{sss:intexplicit} and~\ref{sss:defaut} for an explicit 
  description of the results below in terms of the original datum.
\end{remark}

\subsubsection{Other derived structure}\label{sss:structure}
The following additional sets are obtained from the internality datum: Let 
\begin{equation*}
  \ti{F}=\X\x_D\X
\end{equation*}
Given \((x,y)\in\ti{F}\), the maps \(f_x\) and \(f_y\) have the same image in 
\(\Cc\), so the composition \(\autQ{x}{y}\) is defined.  Thus each point of 
\(\ti{F}\) induces a bijection of \(\Q\) on itself.  We let \(F\) be the 
canonical family for \(\ti{F}\) (i.e., the quotient of \(\ti{F}\) obtained by 
identifying pairs that induce the same map), so that \(F\) is a family of 
bijections of \(\Q\) on itself.

Similarly, any two elements \(x,y\in\X\) give rise to a bijection 
\(\autC{x}{y}\) from \(\Cc_{\pi(x)}\) to \(\Cc_{\pi(y)}\). We denote the 
canonical family for this family by \(H\)\index{aa@\(H\)}:
\begin{equation*}
  H=\X\x\X/E
\end{equation*}
where \((x,y)E(z,w)\) if \(\autC{x}{y}=\autC{z}{w}\).

There are two definable maps \(\pi_d\) and \(\pi_i\) from \(H\) to \(D\), for 
the domain and image of the element:
\begin{gather*}
  \pi_d(\autC{x}{y})=\pi(x)\\
  \pi_i(\autC{x}{y})=\pi(y)
\end{gather*}
We denote the fibres of these maps over \(x\in{}D\) by \(H_x\) and \(H^x\), 
respectively.

Note that if \(\Cc\) is stably embedded, then \(H\), being a family of maps 
between subsets of \(\Cc\), has a canonical injective map into \(\Eq{\Cc}\).  

Given an element \(h\) in \(H_d\) (for some \(d\in{}D\)), we may compose it 
with a function \(f_x\), where \(x\in\X_d\) to obtain a new function 
\(h\circ{}f_x\) from \(\Q\) to \(\Cc\).  We denote the definable set of all 
functions acquired in this way by \(\Wbar{\X}\)\index{aa@\(\Wbar{\X}\)}:
\begin{equation*}
  \Wbar{\X}=H\x_D\X/E
\end{equation*}
where \(E\) is again the relation of defining the same map.  Note that we 
have a canonical map from \(D\) to \(H\) corresponding to the identity map on 
each fibre \(\Cc_a\), hence we have a definable injective map from \(\X\) to 
\(\Wbar{\X}\), and \emph{we consider \(\X\) to be a subset of \(\Wbar{\X}\) 
via this map}. We denote by
\begin{equation*}
  \mu:H\x_D\X\ra\Wbar{\X}
\end{equation*}\index{aa@\(\mu\)}
the quotient map. Note that \(H\), \(\mu\) and \(\Wbar{\X}\) can all be 
defined in the same way in the situation of remark~\ref{rmk:oldint}, if we 
properly modify \(D\), as explained there.

\subsubsection{Definability of the group}\label{sss:defgp}
Following~\ref{sss:images} and~\ref{sss:structure}, the internality datum 
consists of definable sets \(\Q\), \(\X\), \(\Cc\) and \(D\), definable maps 
\(\pi:\X\ra{}D\) and \(\Cc\ra{}D\) onto \(D\), and a definable map 
\(f:\Q\x\X\ra\Cc\) over \(D\), such that the combined map 
\((f,p_2):\Q\x\X\ra\Cc\x_D\X\) (with \(p_2\) the second projection) is 
bijective. We view elements \(x\) of \(\X\) as bijections \(f_x\) from \(\Q\) 
to \(\Cc_{\pi(x)}\), and we assume that distinct elements of \(\X\) given 
distinct maps. We have the auxiliary definable sets \(F\) of bijections from 
\(\Q\) to \(\Q\), and \(H\) of bijections between fibres of \(\Cc\) over 
\(D\), both obtained by composing two elements of \(\X\), and an extra set 
\(\Wbar{\X}\) containing \(\X\) of bijections from \(\Q\) to fibres in 
\(\Cc\), obtained by composing elements of \(H\) with elements of \(\X\), 
with composition map \(\mu:H\x_D\X\ra\Wbar{\X}\).

Given \(x,y\in\X(M)\), we write \(x\Conj{}y\) if there exists an automorphism 
\(\tau\in{}Aut_{(\X,f)}(\Q/\Cc)(M)\) such that \(\tau(x)=y\). We note that 
such \(\tau\) is unique. To see this, it is enough to show that if 
\(\tau(x)=x\) then \(\tau\) is the identity. However, for any \(q\in\Q\),
\begin{equation*}
  f(q,x)=f(\tau(q),\tau(x))=f(\tau(q),x)
\end{equation*}
Since each \(f_x\) is injective, this shows that \(\tau(q)=q\), so \(\tau\) 
is the identity on \(\Q\). This means that for any \(y\in\X\), \(f_y\) and 
\(f_{\tau(y)}\) are the same map on \(\Q\). By the definition of internality, 
this implies that \(\tau(y)=y\), so that \(\tau\) is the identity on \(\X\) 
as well.

\begin{prop}\label{prp:autispair}
  \mbox{}
  \begin{enumerate}[a.]
    \item
      Given an internality datum \((\X,f)\) there is a definable family 
      \(\phi(x,h)\) of subsets of \(H\) parametrised by \(\X\), such that for 
      any model \(M\) and \((x,y)\in\ti{F}(M)\), \(x\Conj{}y\) if and only if 
      they define the same subset of \(H\).

    \item
      If \(\Q\) is \((\X,f)\) internal to \(\Cc\), then the group 
      \(G=Aut_{(\X,f)}(\Q/\Cc)\) is definable: there is a definable group 
      action \(G_0\x\Q\ra\Q\) such that for any model \(M\), \(G(M)=G_0(M)\), 
      with the given action. The definition of the group and the action is 
      given explicitly in terms of the internality datum.

    \item
      In particular, if \(\Q\) is \(\Cc\)-internal relatively to \(G\), then 
      the \(G\) action on \(\Q\) is, in the same sense, a sub-action of a 
      definable group action.
  \end{enumerate}
\end{prop}
\begin{proof}
  \mbox{}
  \begin{enumerate}[a.]
    \item
      A basic observation is that \(\mu\) and \(f\) commute: For any 
      \(q\in\Q\), \(h\in{}H\) and \(x\in\X\), if \(\mu(h,x)\in\X\), then we 
      have
      \begin{equation*}
        f(q,\mu(h,x))=h(f(q,x))
      \end{equation*}
      where we denote by \(h(-)\) the action of \(H\) on \(\Cc\). This is 
      simply the definition of \(\mu\).

      We now claim that the required formula is \(\mu(h,x)\in\X\), i.e., 
      \(x\Conj{}y\) if and only if
      \begin{equation*}
        f_z\circ{}g_w\circ{}f_x\in\X\iff f_z\circ g_w\circ f_y\in\X
      \end{equation*}
      for all \(z,w\in\X\) for which the composition makes sense (where a 
      function ``belongs'' to \(\X\) if it has the form \(f_t\) for some 
      \(t\in\X\)).
      
      In fact, if both \(x\) and \(\mu(h,x)\) belong to \(\X\), and \(\tau\) 
      is an automorphism of the internality structure, then 
      \(\tau(\mu(h,x))=\mu(h,\tau(x))\): to show this, it is enough to show 
      that they coincide as maps from \(\Q\) to \(\Cc\). But
      \begin{equation}
        \begin{split}
      f(q,\tau(\mu(h,x)))&=f(\inv{\tau}(q),\mu(h,x))=h(f(\inv{\tau}(q),x))=\\
      &=h(f(q,\tau(x)))=f(q,\mu(h,\tau(x)))
        \end{split}
      \end{equation}
      Therefore, if \(x\Conj{}y\), then they define the same subset.

      Conversely, assume that \(x\) and \(y\) define the same subset of 
      \(H\). Let \(\tau\) be the map on \(\Q\) given by \(g_y\circ{}f_x\).  
      To show that this map extends to an automorphism of the internality 
      structure, we need to show that given any \(z\in\X\), the map 
      \(f_z\circ{}g_x\circ{}f_y\) from \(\Q\) to \(\Cc\) coincides with 
      \(f_w\) for some (unique) \(w\in\X\). Let \(h\in{}H\) be the element 
      corresponding to \(f_z\circ{}g_x\). Then we need to show that 
      \(\mu(h,y)\in\X\). But \(\mu(h,x)=z\in\X\), so by assumption 
      \(\mu(h,y)\) belongs to \(\X\) as well.

    \item
      As described above, the group is obtained as the subset of \(F\) 
      corresponding to the pairs \((x,y)\in\ti{F}\) such that \(x\Conj{}y\), 
      i.e., pairs that satisfy
      \begin{equation*}
        \forall h\in H(\mu(h,x)\in\X\iff\mu(h,y)\in\X)
      \end{equation*}
      And the action on \(\Q\) is the restriction of the corresponding action 
      of \(F\) (as a set of bijection of \(\Q\) with itself).

    \item
      By definition, \(G\) is a subgroup of the automorphism group of an 
      internality datum for it.\qedhere
  \end{enumerate}
\end{proof}

Thus, if we want to show that the groups \(G(M)\) are the groups of points of 
an \ogroup, we already have a natural map into a definable group, and we only 
need to determine whether the image is a bounded intersection of some of its 
definable subgroups. The answer to this question is given 
in~\ref{ssc:partial}.

\begin{example}\label{ex:groupint}
  The sets \(\Q\), \(\X\) and \(\Cc\) need not be distinct. For example, if 
  \(T\) is the theory of groups, we may take \(\Q=\X=\Cc\) to be the universe 
  \(G\), and \(f:\Q\x\X\ra\Cc\) to be the group multiplication. In this case, 
  the automorphisms group is \(G\) itself, and an element \(g\in{}G\) acts on 
  \(\Q=G\) by \(q\mapsto{}qg^{-1}\) and on \(\X=G\) by \(x\mapsto{}gx\).
\end{example}

We recall the following definition.
\begin{definition}\label{def:torsor}\index{torsor}
  Let \(G\) be a (pro-)definable group. A (pro-)definable group action 
  \(m:G\x{}X\ra{}X\) is called a \emph{\(G\)-torsor} if \(X\) has a point in 
  some model, and the combined map \((p_2,m):G\x{}X\ra{}X\x{}X\) (where 
  \(p_2\) is the projection on the second factor) is an isomorphism
\end{definition}

The claim of proposition~\ref{prp:autispair} can thus be stated as follows: 
The action of \(G\) on \(\X\) is, by definition, free. An orbit of this 
action is given by a definable (over parameters) subset of \(\X\). Any such 
orbit is thus a \(G\)-torsor. If \(\Cc\) is stably embedded, the family of 
all such orbits is in \(\Eq{\Cc}\).

\subsubsection{An explicit definition}\label{sss:intexplicit}
Let us write down the explicit first order definition of the internality 
group.\index{internality group!first order definition of} Following through 
the proof we see that it is (in the free variable \(u\)):
\begin{equation}\label{eq:intdef1}
  \forall(x,y)\in\ti{F}(\Pi(x,y)=u\implies
    \forall h\in H(\mu(h,x)\in\X\iff \mu(h,y)\in\X))
\end{equation}
where \(\Pi\) is the projection from \(\ti{F}\) to \(F\).

Alternatively, we may pass to the canonical family \(F\) first. To do this, 
recall that \(F\) is a family of bijections from \(\Q\) to itself, and let 
\(\what{\X}\)  be the canonical family of bijections from \(\Q\) into \(\Cc\) 
obtained by composing elements of \(F\) and \(\X\). Let
\begin{equation*}
  \nu:F\x\X\ra\what{\X}
\end{equation*}
be the composition map, and, as for \(\bar{\X}\), we identify \(\X\) with the 
subset \(\nu(id,\X)\) of \(\what{\X}\). More generally, we identify elements 
of \(\what{\X}\) and of \(\Wbar{\X}\) that define the same function from 
\(\Q\) to \(\Cc\). After this identification, if \(h\in{}H\), \(u\in{}F\) and 
\(x\in\X\) are elements with the property that \(\mu(h,x),\nu(u,x)\in\X\), 
then \(\nu(u,\mu(h,x))=\mu(h,\nu(u,x))\), since they \(\mu\) and \(\nu\) 
correspond to composition on different sides.  We also note that given any 
\(u\in{}F\), the inverse of \(u\) (as a function on \(\Q\)) is also 
represented by an element of \(F\), denoted \(\inv{u}\).
\begin{claim}\label{clm:altdef}
  The internality group given by the subset of elements \(u\in{}F\) 
  satisfying the formula
  \begin{equation}\label{eq:intdef2}
    \forall z\in\X(\nu(u,z)\in\X\land\nu(\inv{u},z)\in\X)
  \end{equation}
\end{claim}
\begin{proof}
  Assume that \(u=\Pi(x,y)\) satisfies formula~\eqref{eq:intdef1}. Given 
  \(z\in\X\), let \(h\in{}H\) be the element corresponding to 
  \(\autC{x}{z}\).  Then \(\mu(h,x)=z\) and \(\mu(h,y)=\nu(u,z)\). Hence, by 
  formula~\eqref{eq:intdef1}, \(\nu(u,z)\in\X\). The proof for the inverse is 
  symmetric.

  Conversely, if \(u\) satisfies formula~\eqref{eq:intdef2} and 
  \(\Pi(x,y)=u\), assume that \(\mu(h,x)\in\X\) for some \(h\in{}H\). Then 
  for \(z=\mu(h,x)\), we get:
  \begin{equation*}
    \mu(h,y)=\mu(h,\nu(u,x))=\nu(u,\mu(h,x))=\nu(u,z)\in\X
  \end{equation*}
  and the other direction follows from the condition on \(\inv{u}\).
\end{proof}

We would like also to rewrite this last formula explicitly in terms of the 
original data. Thus, instead of \(u\in{}F\), we describe its pre-image in 
\(\X\x\X\) in terms of the function \(f\) (which need not be given via a 
function symbol, but as a ternary relation.) Expanding the definition we get 
(in the free variables \((z,w)\in\X\x\X\)):
\begin{gather}\label{eqn:prehorrible}
  \forall x\in\X\exists y\in\X\phi(x,y,z,w)\\
  \intertext{where the formula \(\phi(x,y,z,w)\), representing the condition  
  \(f_y=\inv{f_w}\circ\autC{x}{z}\), is given by}
  \begin{split}\label{eqn:horrible}
    \forall q\in\Q,c\in\Cc( & f(q,y,c)\iff\\
    &\exists p\in\Q,d\in\Cc(f(q,w,d)\land f(p,z,d)\land f(p,x,c)))
  \end{split}
\end{gather}
This equation is valid also when \(f\) has the form in 
remark~\ref{rmk:oldint}. In these terms, the action of \(G\) on \(\Q\) is 
given (via \(\X\x\X\)) by \((z,w)(p)=q\), where \(z,w\in\X\) and \(p,q\in\Q\) 
satisfy
\begin{equation}\label{eqn:horribleQ}
  \exists c\in\Cc(f(q,w,c)\land f(p,z,c))
\end{equation}
and the action on \(\X\) is given by \((z,w)(x)=y\), where \(x,y,z,w\in\X\) 
satisfy \(\phi\) of equation~\eqref{eqn:horrible}.

We stress that these definitions are valid in any model. Thus, the above 
formulae describe the group of automorphisms over \(\Cc\) of the datum  
\(f:\Q\x{}\X\ra{}\Cc\) where \(\Q\), \(\X\) and \(\Cc\) are arbitrary sets, 
and \(f\) is an arbitrary function satisfying the axioms of internality.

\subsubsection{A minimal example}\label{sss:intexample}
\index{internality!minimal example of}
The following example deals with a minimal situation. The language contains 
three sorts, \(\Q\), \(\X\) and \(\Cc\), and two function symbols, 
\(f:\Q\x\X\ra\Cc\) and \(g:\Cc\x\X\ra\Q\). The theory \(T\) says that the 
maps \((f,p_2):\Q\x\X\ra\Cc\x\X\) and \((g,p_2):\Cc\x\X\ra\Q\x\X\) are 
bijective and inverse to each other (as before, \(p_2\) is the projection to 
the second component.) This is a universal theory, and the theories 
considered below extend \(T\). For each of them we will consider the 
following possibilities for \(G\): \(G_0(M)\) is the group of all bijections 
of \(\Q(M)\) preserving the quantifier free subsets of all \(\Q^n\) definable 
with parameters from \(\Cc(M)\). \(G_1(M)\) is the group preserving all 
\(\Cc(M)\) definable subsets of all \(\Q^n\), and \(G_2(M)\) is the full 
group of automorphism of \(M\) over \(\Cc(M)\), restricted to \(\Q\): 
\(G_2(M)=Aut(M/\Cc(M))/Aut(M/\Q(M),\Cc(M))\).

The group \(G_0\) can be described immediately: \(\Q\) has no quantifier free 
structure at all over any set of parameters contained in \(\Cc\). Therefore, 
\(G_0(M)\) is simply the set of bijections of \(\Q(M)\) onto itself as a set.

\begin{enumerate}
  \item
    Let \(T_1\) be the theory saying that for any two elements 
    \((q,c)\in\Q\x\Cc\), there is a unique element of \(\X\) mapping \(q\) 
    to \(c\).  Clearly, most elements of \(G_0(M)\) can not be extended to 
    \(\X\).

    Let \(\tau\in{}G_1(M)\), and let \(q_0\in\Q(M)\) be any point. Extend 
    \(\tau\) to \(\X(M)\) by setting 
    \(f(q_0,{}\tau(x))=f(\inv{\tau}(q_0),{}x)\).  To show that this is an 
    automorphism, we need to show that the same holds for any other point 
    \(q\in\Q(M)\). By the axioms, any element \(c\in\Cc(M)\) gives rise to a 
    definable bijection between \(\Q\) and \(\X\) (namely \(g(c,-)\).) Given 
    two such elements, we get a bijection of \(\Q\) to itself, definable over 
    these two elements. In particular, this bijection commutes with \(\tau\). 
    Considering this bijection for the two elements \(x(q)\) and \(x(q_0)\) 
    we get the result. Thus, this is an internality datum for \(G_1\). Also, 
    since the whole structure consists of the maps \(f\) and \(g\), we see 
    that \(G_1\) and \(G_2\) coincide.

    One may also consider the theories \(T_n\), which say that for any 
    pairwise distinct \(q_1,\dots,q_n\in\Q\) and \(c_1,\dots,c_n\in\Cc\) 
    there is a unique \(x\in\X\) mapping \(q_i\) to \(c_i\). The same 
    description as for \(T_1\) holds there. In fact, any model of \(T_n\) 
    gives rise to a model of \(T_1\) with the same groups by considering the 
    subsets of \(\Q^n\) and \(\Cc^n\) consisting of tuples with pairwise 
    distinct coordinates.

  \item
    Let \(M_0\) be the model of \(T\) where \(Q_0\) and \(C_0\) are infinite 
    sets of the same cardinality, and \(X_0\) the set of all bijections from 
    \(Q_0\) to \(C_0\). Let \(T_{\Inf}\) be the theory of \(M_0\).  For this 
    particular model \(M_0\), \(G_0(M_0)=G_1(M_0)=G_2(M_0)\) all coincide, 
    since the action of the group \(G_0(M_0)\) of all bijection of 
    \(\Q(M_0)\) to itself extends uniquely to an action on \(\X(M_0)\) (by 
    composition), so that and \(G_0(M_0)=G_2(M_0)\), and \(G_1\) is between 
    the two groups.
    
    In particular, since \(G_2(M_0)\) acts transitively on the set of 
    pairwise distinct \(n\)-tuples (for any \(n\)), the \(\Cc(M)\) 
    definable subsets of \(\Q\) are quantifier free, for any model \(M\) of 
    \(T_{\Inf}\).  Therefore, \(G_0\) and \(G_1\) coincide completely for 
    this theory. 
    
    However, not every model is of the same form as \(M_0\): by the 
    omitting types theorem, we may find a model \(M_1\) with \(Q_1=Q_0\), 
    \(C_1=C_0\) but where \(X_1\) is a proper subset of \(X_0\). Therefore 
    \(G_2(M_1)\) is strictly smaller than \(G_1(M_1)\). This shows that the 
    condition that automorphisms should be extendible to \(\X\) is 
    non-trivial, even when there is quantifier elimination for subsets of 
    \(\Q\). There are some more remarks about this situation in 
    section~\ref{ssc:partial}.
\end{enumerate}

\subsubsection{Ind-definable internal sets}\label{sss:Qind}
We conclude this sub-section with a remark (which will not be used anywhere) 
about the case when \(\Q\) is ind-definable (See, e.g., \Cite{prodef} for the 
notion. An example is an infinite union of definable sets; an example of a 
pro-definable set is a partial type).\index{internality!for ind-definable 
sets|(} In this case the statement of proposition~\ref{prp:autispair} is 
false, and most of the structure defined here can no longer be defined (there 
are no canonical families of ind-definable sets.) However, if we require the 
definable group to agree with the (slightly modified) group of automorphisms 
for saturated enough models, and use the horrible 
formula~\eqref{eqn:horrible}, we obtain a description.

So we assume now that \(\Q=\Ind[i]{\Q_i}\) and \(\Cc=\Ind[i]{\Cc_i}\) are 
ind-definable sets, that can be presented via systems \(\Q_i\) and \(\Cc_i\) 
of size less than a cardinal \(\kappa\). \(\X\) is still a definable set, and 
we are given an ind-definable map \(f:\Q\x\X\ra\Cc\), which amounts to giving 
a system of definable maps \(f_i:\Q_i\x\X\ra\Cc_j\) (for every \(i\) for some 
\(j\).) However, the requirement that the \(f_x\) are distinct is no longer 
first order. We therefore modify the definition of 
\(Aut_{(\X,f)}(\Q/\Cc)(M)\) to be the set of pairs \((\tau_{\Q},\tau_{\X})\), 
where \(\tau_{\Q}\), as before, is a (set) automorphism of \(\Q(M)\), but 
\(\tau_{\X}\) is a bijection of \(\Wbar{\X(M)}\), where \(\Wbar{\X(M)}\) is 
the set of maps from \(\Q(M)\) to \(\Cc(M)\) represented by elements of 
\(\X(M)\) (we could have taken the same path when defining internality in the 
first order case, but the existence of \(\EQ{T}\) allows us to use the 
simpler form.) We assume that the maps in the system \(\Q_i\) are injective.  
It follows that the injectivity of each fibre \(f_x\) is a first order 
property (we could avoid this assumption by further modifying the definition 
of the automorphisms group into something unrecognisable, but for simplicity 
we leave it as it is.)

We note that the properties described in the formulas~\eqref{eqn:horrible} 
and~\eqref{eqn:horribleQ} still define the group and the actions, but they 
are no longer first order. Rewriting them in a way that the quantifiers 
become first order, we get instead of \(\phi\) there:
\begin{gather}\label{eqn:veryhorrible}
  \Pro[i]{\forall q\in\Q_i,c\in\Cc_{j_i}(f_i(q,y,c)\iff
  \Ind[k]{\psi_{i,k}(z,w,q,x,c)})}\\
  \intertext{where \(\psi_{i,k}(z,w,q,x,c)\) is}
   \exists p\in\Q_k,d\in\Cc_{l_k}(f_i(q,w,d)\land f_k(p,z,d)\land f_k(p,x,c))
\end{gather}
At this point it is convenient to split the double implication (\(\iff\)) 
into its two components (with which the universal quantifiers on \(q\) and 
\(c\) commute.) After doing this and rearranging terms, we get that a typical 
term in the projective system above is the intersection of the two 
expressions:
\begin{gather}
  \forall q,c(\Ind{f(q,y,c)\implies\psi})\\
  \forall q,c(\Pro{\psi\implies f(q,y,c)}
\end{gather}
The universal quantifiers always commute with projective systems, but they 
commute with the inductive ones only for \(\kappa\)-saturated models. For 
such models we thus get an intersection of a pro-definable set with an 
ind-definable one (all subsets of \(\X^4\)), which can be viewed as a 
pro-ind-definable. Plugging this into the formula~\eqref{eqn:veryhorrible}, 
we get again a pro-ind-definable formula (which is the formula for the action 
of the automorphisms group on \(\X\).) Finally, the \(\kappa\)-saturation 
allows us again to move the quantifiers on \(x\) and \(y\) in 
formula~\eqref{eqn:prehorrible} inside, resulting in a pro-ind-definable 
definition for the set of pairs of elements of \(\X\) representing 
automorphisms. The automorphisms group is now obtained by dividing by a 
projective system of equivalence relations, so the final answer is that the 
automorphisms group coincides, for \(\kappa\)-saturated models, with a 
pro-ind-definable group (this is not very surprising --- we knew in advance 
that any individual automorphism is definable, so the group is at least a 
subset of the pro-ind-definable set \(Hom(\Q,\Q)\) of all definable maps from 
\(\Q\) to \(\Q\).)

As an example, consider the theory, in the language of internality 
(\ref{sss:intexample}), of the model where \(\Cc\) and \(\X\) are the 
rational numbers, and \(\Q\) is the ind-definable set of ``finite'' rational 
numbers, the union of intervals \((-n,n)\), and \(f\) is given by the 
restriction of addition. Then the formula~\eqref{eqn:veryhorrible} defines 
the set of pairs \((z,w)\) of elements of \(\X\) whose distance is finite, 
and the group is the ind-definable group of finite numbers. This example 
shows, in particular, that although each \(f_i\) gives rise to an internality 
datum on \(\Q_i\), the resulting groups have nothing in common with the group 
of the whole system.

On the other hand, if the system \(f\) has the property that for any \(j\) 
there is an \(i\) such that if \(f(q,x)=c\) for some \(x\in\X\) and 
\(c\in\Cc_j\), then \(q\in\Q_i\), then the group is pro-definable. If, 
furthermore, the restriction of \(f\) to \(\Cc_i\) is just \(f_i\), then the 
resulting group is an \ogroup, whose system is given by the groups 
corresponding to the \(f_i\).\index{internality!for ind-definable sets|)}
\index{internality|)}

\subsection{Partial automorphisms}\label{ssc:partial}
We are going to be concerned with a class of definable sets, \(\QF\), and 
restate some of the definitions relative to this class. \(\QF\) sets can be 
thought to be quantifier free, in the sense that they \emph{need not} be 
closed under applying quantifiers (and also because this class is a central 
example.) However, no actual properties of quantifier free sets will be used.  
The collection \(\QF\) \emph{will} be assumed to be closed under boolean 
combinations, but only to make the formulation simpler.
\begin{definition}
  Let \(S=\{X_i\}\) be any collection of definable sets, \(M\) a model. An 
  \emph{automorphism}\index{automorphism!of a collection of definable sets} 
  of \(S(M)\) is a collection of invertible maps
  \begin{equation*}
    f_i:X_i(M)\ra X_i(M)
  \end{equation*}
  such that for any Cartesian product \(Y\) of elements of \(S\), the induced 
  map preserves all \(\QF\) subsets of \(Y\).
\end{definition}
We note that equality is not required to be in \(\QF\), so, as in 
example~\ref{ex:groupint}, a definable set may appear more than once in the 
list \(X_i\), with different \(f_i\).

We denote the group of all such automorphisms\index{automorphisms group} by 
\(\Aut{S}(M)\). When \(T\subseteq{}S\), we have a natural restriction map 
from \(\Aut{S}(M)\) to the group \(\Aut{T}(M)\), and we will be interested in 
the kernel \(\Aut[T]{S}(M)\)\index{aa@\(\Aut[T]{S}(M)\)} of this map (i.e., 
the automorphisms of \(S\) that preserve \(T\) pointwise.)

Thus, if the set \(S\) consists of all the sorts, and \(\QF\) contains all 
quantifier free sets, then this is what is usually called \(Aut(M)\): The 
group of all automorphisms of \(M\). In contrast, we are interested in 
automorphisms of some of the sorts, and with respect to part of the 
structure.

The following examples show that in general, an automorphism of one sort does 
not extend to other sorts, even if it preserves all of the quantifier free 
structure.

\begin{example}\label{xpl:noextend}
  Let \(T\) be the theory of groups (in the natural language) with an extra 
  predicate \(X\) for a subgroup of index \(2\). Consider the group 
  \(M=\Z\x2\Z\), and let \(S\) contain only the sub-group 
  \(X(M)=2\Z\x2\Z\). Then the function that swaps the coordinates is an 
  automorphism of the quantifier free structure on \(X\) (which is just the 
  group structure), but not of the full structure, since it does not 
  preserve the set \(\exists{}x(y=x+x)\).
\end{example}

\begin{example}
  Let \(T\) be a theory saying that a sort \(k\) is an algebraically closed 
  field, and \(U\) is another sort, such that \(U\x{}U\) is a finite 
  dimensional vector space over \(k\), and the function that swaps the 
  coordinates of \(U\x{}U\) is a linear transformation. If \(\QF\) consists 
  of just the linear structure (in particular, it does not contain the 
  projections to \(U\)), then the swap of coordinates is an automorphism of 
  \(U\x{}U\) (in the sense defined above), that does not extend to an 
  automorphism of \(U\) (if \(U\) has more than one element.)
\end{example}

Additional examples are provided by the theories in \ref{sss:intexample}, 
as well as the example of \ACFA{}, studied in more details in 
section~\ref{sec:ACFA}.

\subsubsection{Relation with the full automorphisms group}\label{sss:compare}
If \(\QF\) is the set of all definable sets, \(X\) is stably 
embedded\index{stably embedded} and \(M\) is saturated, then any automorphism 
\(f\) of \(X(M)\) can be extended to \(M\): Indeed, first we may assume that 
\(T\) eliminates imaginaries, since \(f\) extends uniquely to \(\EQ{X}\).  
Next, assume we managed to extend \(f\) to some small subset 
\(A{\subseteq}{}M\), and we want to extend it further to \(a\).  Let 
\(p=\Tp{a}{X(M)\cup{}A}\).  Then \(f(p)\) is consistent, and as explained 
in~\ref{ssc:SE}, determined by its restriction to a small subset of \(X(M)\) 
and \(A\). By saturation it is realised in \(M\) by some \(b\), and we may 
set \(f(a)=b\).

In particular, under the above assumptions,
\begin{equation*}
  G=\Aut[X(M)]{M}
\end{equation*}
acts transitively on the realisations in \(M\) of any type over \(X(M)\), and 
the definable closure \(dcl(X(M))\) is equal to \(M^G\), the set of fixed 
points of the \(G\) action on \(M\). This is another instance of the 
similarity between stably embedded sets and small sets (The converse of this 
fact is also true: if \(X\) is a definable set such that any automorphism of 
\(X\) lifts to an automorphism of \(M\), then \(X\) is stably embedded.  
Cf.~\Cite[appendix]{ACFA}.)

In this context, it is fruitful to reconsider the example 
in~\ref{sss:intexample}\index{internality!minimal example of}: For \(T_1\), 
the graph of an element \(x\in\X\) as a function from \(\Q\) to \(\Cc\) is 
definable using parameters from \(\Q\) and \(\Cc\) (namely, an arbitrary 
point on that graph.) In fact, \(\Q\cup\Cc\) is stably embedded in this 
theory.  Therefore, \(G_1\) and \(G_2\) coincide.  On the other hand, in 
\(T_{\Inf}\) the graph of an element of \(\X\) is not definable using 
parameters from \(\Q\) and \(\Cc\), and \(G_1\) is bigger than \(G_2\). A 
saturated model is indeed not of the form where \(\X\) is the set of all 
bijections from \(\Q\) to \(\Cc\) (since the set of such bijections for a 
model \(M\) has cardinality strictly bigger than the cardinality of 
\(\Q[M]\); hence given \(x_0\in\X[M]\), its type \(p(x)\) over 
\(\Q[M]\cup\Cc[M]\), which describes the bijection completely, would be of 
cardinality smaller than the cardinality of \(M\), so the type 
\(p(x)\cup\{x\ne{}x_0\}\) would have to be realised in \(M\), and this is a 
contradiction.)

\subsubsection{Definability of the automorphisms group}\label{sss:defaut}
A \emph{\(\QF\)-type}\index{type@\(\QF\)-type} over a set \(A\) is the 
restriction of a usual type over \(A\) to formulas in \(\QF\).

We may now answer the question raised in section~\ref{ssc:Internality}: Which 
subgroups of the internality automorphism group are \ogroups? The following 
is an analogue of proposition~\ref{prp:autispair}.
\begin{prop}\label{prp:defaut}
  \mbox{}
  \begin{enumerate}[a.]
    \item
      Let \(\Q\) be \((\X,f)\)-internal to an (ind-)definable set \(\Cc\), 
      \(\QF\) a set of formulas in the sorts \(\Q\), \(\X\) and \(\Cc\), with 
      \(f\in\QF\), and let \(G=\Aut[\Cc]{\Q,\X}\).  Then there is a set of 
      formulas \(\QF^*\) with one \(\X\) variable and no \(\Q\) variables, 
      such that the \(G\) orbit of an element in \(\X[M]\) is a 
      \(\QF^*\)-type over \(\Cc[M]\cup{}H(M)\cup{}D(M)\).  It thus follows 
      that \(G\)\index{automorphisms group} is an 
      \ogroup.\index{group@\ogroup}

    \item
      Let \(\Q\) be a definable set, \(\Cc\) an (ind-)definable set, \(\QF\) 
      a collection of formulas, such that \(\Q\) is internal to \(\Cc\) 
      relatively to \(G=\Aut[\Cc]{\Q}\). Then \(G\) is an \ogroup.
  \end{enumerate}
\end{prop}

In the proof we shall use the internality datum, and the derived functions 
\(g\) and \(\mu\), to convert formulas
\begin{equation}\label{eqn:phi}
  \phi(q_1,\dots,q_n,x_1,\dots,x_l,\Wbar{e})
\end{equation}
with variables \(\Wbar{q}\in\Q\), \(\Wbar{x}\in\X\) and \(\Wbar{e}\in\Cc\) 
into formulas in \emph{one} \(\X\) variable \(x\) and tuples 
\(\Wbar{d}\in\Cc\) and \(\Wbar{h}\in{}H\), as well as the original 
\(\Wbar{e}\). Given a formula as above, \(\phi^*\) is the subset of
\begin{equation*}
  \X\x_D\Cc\x_D\dots\x_D\Cc\x_D H\dots\x_D H\x\dots
\end{equation*}
given by
\begin{equation}\label{eqn:phistar}
  \phi(g(d_1,x),\dots,g(d_n,x),\mu(h_1,x),\dots,\mu(h_l,x),\Wbar{e})
\end{equation}

\begin{proof}[Proof of~\ref{prp:defaut}]
  Note that given an internality datum \((\X,f)\) and an action of \(G\) on 
  \(\X\), saying that \(f\) belongs to \(\QF\) is equivalent to saying that 
  \((\X,f)\) is an internality datum for \(G\). Therefore, the second part is 
  a corollary of the first.

  Let \(\QF_1=\QF\cup\{\X\}\) (\(\X\) is the formula \(x\in\X\)). Any formula 
  \(\phi\) in \(\QF_1\) is of the form~\eqref{eqn:phi}, and we set
  \(\QF^*=\{\phi^*\|\phi\in\QF_1\}\), where \(\phi^*\) is as 
  in~\eqref{eqn:phistar}. Note that, in particular, \(\X^*\in\QF^*\) is the 
  formula \(\mu(h,x)\in\X\), which is the formula \(\phi\) promised by 
  proposition~\ref{prp:autispair}.

  We need to show that given \((z,y)\in\ti{F}\), \(z\) and \(y\) have the same 
  \(\QF^*\)-type if and only if they are in the same orbit of \(G\). Since 
  \(\X\in\QF_1\), they are in the same orbit of the internality group. Let 
  \(\tau\) be the element in that group taking \(z\) to \(y\). Let 
  \(\phi(\Wbar{q},\Wbar{x},\dots)\in\QF\).  We claim, similarly to 
  proposition~\ref{prp:autispair}, that \(\tau\) preserves \(\phi\) if and 
  only if \(\phi^*(y,\Wbar{c})\) and \(\phi^*(z,\Wbar{c})\) are the same set 
  (where \(\Wbar{c}\) contains everything except the variable in \(\X\).) The 
  proof is also similar:

  Let \(\Wbar{q}\in\Q[M]\), \(\Wbar{x}\in\X[M]\) be tuples of elements, let  
  \(\Wbar{d_1}\) be the elements of \(H\) obtained by the compositions 
  \(f_{x_i}\circ{}g_z\), and let 
  \(\Wbar{c_1}=f(\Wbar{q},z)=(f(q_1,z),\dots,f(q_n,z))\). Since \(\tau\) is 
  an automorphism of the internality structure, we also have 
  \(\Wbar{c_1}=f(\tau(\Wbar{q}),y)\) and \(\Wbar{d_1}\) is also the element 
  corresponding to \(f_{\tau(\Wbar{x})}\circ{}g_y\). On the other hand, by 
  the definition of \(g\) we have \(g(\Wbar{c_1},z)=\Wbar{q}\) and 
  \(g(\Wbar{c_1},y)=\tau(\Wbar{q})\).  Likewise, we have 
  \(\mu(\Wbar{d_1},z)=\Wbar{x}\) and \(\mu(\Wbar{d_1},y)=\tau(\Wbar{x})\).  
  Thus, if \(z\) and \(y\) define the same \(\phi^*\) subset, then \(\tau\) 
  preserves \(\phi\). The converse is also true, since by definition, every 
  element of \(\Cc_{\pi(z)}\) is the image of some \(q\in\Q\) under the 
  action of \(z\), and for every \(x\in\X\) we may take \(d\in{}H\) 
  corresponding to \(f_x\circ{}g_z\).
\end{proof}

Similarly to the situation of pure internality (\ref{sss:intexplicit}), the 
definition of the group is explicit in terms of \(\QF\) and the internality 
datum. Let \(\Pi\) be the projection from \(\ti{F}\) to \(F\).  Then the 
definition of the automorphisms group
\index{automorphisms group!first order definition of} is given explicitly 
(in the free variable \(u\)) as
\begin{equation*}
\bigand{\phi\in\QF^*}\forall(x,y)\in\ti{F}
(\Pi(x,y)=u\implies\forall\Wbar{c}(\phi(x,\Wbar{c})\iff\phi(y,\Wbar{c})))
\end{equation*}
where \(\QF^*=\{\phi^*\|\phi\in\QF_1\}\) is, as before, the set of \(\QF\) 
formulas composed with \(g\) and \(\mu\), as given by 
equation~\ref{eqn:phistar}.
In particular, this is a universal formula relatively to the formulas in 
\(\QF^*\), the maps \(\pi:\X\ra{}D\), \(\Pi:\ti{F}\ra{}F\), and the maps 
from \(\Cc\) to \(D\). Note also, that we may use an existential 
quantifier, instead of a universal one:
\begin{equation*}
\bigand{\phi\in\QF^*}\exists(x,y)\in\ti{F}
(\Pi(x,y)=g\land\forall\Wbar{c}(\phi(x,\Wbar{c})\iff\phi(y,\Wbar{c})))
\end{equation*}
This is because the property of having the same \(\QF\) type over \(\Cc\) is 
constant on fibres of \(\Pi\), so that one pair has it if and only if any 
pair in the fibre has it.

Alternatively, as in~\ref{sss:intexplicit}, we may first pass to the 
quotient, and describe the group via its action on \(\X\). We thus get:
\begin{equation*}
u\in F\land
\bigand{\phi\in\QF^*}\forall x\in\X,\Wbar{c}\in\Cc
(\phi(x,\Wbar{c})\iff \phi(\nu(u,x),\Wbar{c}))
\end{equation*}
where \(\nu\) is the function from~\ref{sss:intexplicit}, and whose 
restriction to \(G\subseteq{}F\) is just the action of \(G\) on \(\X\).

We note that, in addition to the quantifiers appearing explicitly in the 
definition, quantifiers can appear in these formulas resulting from the 
transformation of formulas from \(\QF\) to \(\QF^*\), even if the formulas in 
\(\QF\) are quantifier free.

In section~\ref{sec:stability} we obtain, for a particular class of examples, 
a more explicit description, in terms of ``rational function'' from \(\X\) to 
\(\Cc\).

\begin{remark}
  As in the case of pure internality, this result can be restated by saying 
  that the type provided is a torsor\index{torsor} over \(G\). In fact, it is 
  an \(\w\)-torsor\index{torsor@\(\w\)-torsor}, in the sense that it is the 
  intersection of \(G_i\)-torsors, for the groups \(G_i\) preserving finite 
  subsets of \(\QF\), whose intersection is \(G\).
\end{remark}

\subsubsection{Relation with the classical group}\label{srk:original}
We would like to compare the group just constructed with the classical result 
on the binding group. We shall use the setting of appendix B of~\Cite{DG}. 
There, one works with a fixed saturated model \(M\) of a theory with 
quantifier elimination and EI, and one considers the group 
\(\hat{G}(M)=Aut(M/\Cc[M])/Aut(M/\Q[M],\Cc[M])\).  It is shown there that for 
a saturated model \(M\), this group is the group of \(M\) points of an 
\(\w\)-definable group (\(\Q\) is not assumed to be stably embedded.) Note 
that in general, when \(M\) is not saturated, \(\hat{G}(M)\) need not 
correspond to  the \(M\)-points of this \(\w\)-definable groups (in fact, for 
some models \(\hat{G}(M)\) may be trivial.)

What is the connection between this group and the group \(G\) defined here, 
when \(\QF\) is the set of all formulas?

If \(\Cc\) is stably embedded, these are the same groups. Indeed, \(\hat{G}\) 
is obviously contained in \(G\), but as noticed before, \(\hat{G}\) acts 
transitively on the type of an element \(x\in\X\) over \(\Cc\), which we just 
saw to be a \(G\)-torsor. In particular, any automorphism of the full 
structure on \(\Q\) and \(\X\) fixing \(\Cc\) can be extended to the whole 
(saturated) model. This is not, in general, true for automorphisms not fixing 
\(\Cc\) --- the triple \(\Q,\X,\Cc\) need not be stably embedded.

For general \(\Cc\), let \(\Se{\Cc}\) (the \Def{stably embedded hull} of 
\(\Cc\)), be the collection of definable sets \(Y\) (in the whole theory) 
such that  \(Y(M)\) is fixed pointwise by \(Aut(M/\Cc[M])\) (this does not 
actually depend on \(M\)).  Obviously, the automorphisms group \(\hat{G}\) 
does not change when we replace \(\Cc\) by \(\Se{\Cc}\) (this is not true for 
\(G\).)  However, \(\Se{\Cc}\) is stably embedded: given a canonical family 
of subsets of \(\Se{\Cc}\), any automorphism fixing \(\Se{\Cc}\) pointwise 
must fix the parameter of the family as well, so this parameter belongs to a 
set in \(\Se{\Cc}\).

Therefore, in any case \(\hat{G}\) coincides with a group of the form 
considered here.

\subsubsection{Deriving internality data from \(\QF\)}
In the example at the beginning of this section, where \(\Q\) is a 
finite-dimensional vector space over the field \(\Cc\), the internality 
datum, i.e., the set of bases \(\X\), and the linear combinations function 
\(f\) are derived from the linear structure \(\QF\) that we would like to 
preserve. However, with our definitions this need not be the case in general.  
In other words, assuming that \(\Q\) is 
\((\X,f)\)-internal\index{internal!\((\X,f)\)-} to \(\Cc\), and given a set 
of formulas \(\QF\), it need not a-priori be the case that \((\X,f)\) is an 
internality datum for \(G=\Aut[\Cc]{\Q}\).
  
Let \(\LL_\QF\) be the language with the sorts \(\Q,\Cc\), and whose basic 
relations are all definable subsets of \(\Cc\) and all sets in \(\QF\). Let 
\(T_\QF\) be the theory \(T\) restricted to \(\LL_\QF\).  Clearly, any 
definable set in \(\EQ{T_\QF}\) can be identified with a definable set in 
\(\Eq{\Q,\Cc}\).  If \(\X\) and \(f\) correspond in this way to definable 
sets in \(\EQ{T_\QF}\), then they form an internality datum for \(G\), since 
\(G\) is the full automorphism group of \(\EQ{T_\QF}\). Conversely, we have:

\begin{prop} \label{prp:intinEQ}
  Assume that in \(T\) there exists an internality datum for 
  \(G=\Aut[\Cc]{\Q}\), and that \(\Cc\) is stably embedded in \(T_\QF\). Then 
  there exists such a datum for \(G\) in \(\EQ{T_\QF}\).
\end{prop}
\begin{proof}
  Let \(M\) be a model of \(T\), and \(M_\QF\) its restriction to \(T_\QF\).  
  Then \(G(M)\) is the full automorphism group of \(M_\QF\) over 
  \(C(M_\QF)\). Since there is an internality datum for \(G\), \(G\) is an 
  \ogroup (in \(T\).) Therefore, the cardinality of \(G(M)\) is at most that 
  of \(M\) and \(M_\QF\). The result now follows from 
  proposition~\ref{prp:card} below.
\end{proof}

\begin{prop}\label{prp:card}
  Let \(T\) be a theory, \(\Cc\) a definable set stably embedded in \(T\).  
  If there is a saturated model \(M\) of \(T\) such that 
  \(G(M)=Aut(M/\Cc[M])\) has at most the cardinality of \(M\), then the 
  universe is internal to \(\Cc\).
\end{prop}
\begin{proof}
  If there is a tuple \(a\in{}M\), such that any automorphism fixing 
  \(\Cc[M]\) and \(a\) is the identity, then, since \(\Cc\) is stably 
  embedded and \(M\) is saturated, \(dcl(a\cup\Cc[M])=M\). Therefore, \(a\) 
  defines a surjection from a power of \(\Cc\) onto \(M\) (more precisely, 
  any element of \(M\) is in the image of some \(a\) definable map from 
  \(\Cc\). By compactness, a finite number of these maps suffices, and this 
  finite number can be combined into one surjective map on a quotient of the 
  union of their domains, a set in \(\Eq{\Cc}\).)

  We assume there is no such tuple, and will show that the cardinality of the 
  group is bigger than the cardinality \(\ka\) of \(M\). By saturation, it 
  follows that for any subset \(A\) of \(M\) of cardinality less than 
  \(\ka\), and any automorphism \(f\) over \(\Cc\) one may find a different 
  automorphism over \(\Cc\), which agrees with \(f\) on \(A\). This allows us 
  to build a binary tree of height \(\ka\), with different branches 
  corresponding to different automorphisms over \(\Cc\). The number of such 
  branches is \(2^\ka\), so we get a contradiction.
\end{proof}

Considering the example of~\ref{sss:intexample}\index{internality!minimal 
example of} again, we see that indeed in \(T_1\), \(\X\) is a quotient of 
\(\Q\x\Cc\), whereas in \(T_{\Inf}\) (where \(\Cc\) is not stably embedded), 
\(\X\) might have larger cardinality than \(\Q\) and \(\Cc\) in some models, 
hence can not belong to \(T_\QF\).

\subsubsection{A partial converse}
Proposition~\ref{prp:defaut} has the following partial converse:
\begin{prop}\label{prp:defisaut}
  Let \(\Q\) be internal to \(\Cc\), and let \(G\) be an \ogroup acting 
  faithfully on \(\Q\). Then \(\Q\) is \((\X,f)\)-internal to \(\Cc\), where 
  the internality datum is compatible with the action of \(G\), and there is 
  an ind-definable set \(\ti\Cc\) containing \(\Cc\), and a set \(\QF\) of 
  definable subsets of (Cartesian products of) \(\Q,\X,\ti\Cc\), such that 
  \(G=\Aut[\ti\Cc]{\Q,\X}\)
\end{prop}
\begin{proof}
  Let \(G\) be the intersection of a decreasing chain of definable groups 
  \(G_i\). We first claim that the action of \(G\) on \(\Q\) is the 
  restriction of an action on \(\Q\) of some \(G_0\). In fact, by 
  compactness, the action is the restriction of some function 
  \(f:G_0\x\Q\ra\Q\). For any \(i\geq{}0\), let \(\Q_i\) be the set defined 
  by \(\forall{}g,h\in{}G_i(f(g,f(h,q))=f(gh,q))\).  Then the union of all 
  \(\Q_i\) is \(\Q\), a definable set. The same is true of some finite union, 
  hence \(f\) is an action when restricted to some \(G_i\), which from now on 
  is denoted \(G_0\).

  We may now construct \((\X,f)\). Let \((X_1,f_1)\) be the given internality 
  datum, let \(X_2=X_1\x{}G_0\), and let \(f_2:\Q\x{}X_2\MapsTo{}\Cc\) be 
  given by \(f_2(q,x,g)=f_1(gq,x)\). With the action of \(G_0\) on \(X_2\) 
  given by \(h(x,g)=(x,g\inv{h})\), we obtain, by taking the canonical 
  family, internality datum \((\X,f)\) compatible with \(G_0\). Since \(G\) 
  is contained in \(G_0\) this is also an internality datum for \(G\) (and 
  any other \(G_i\).)
  
  Let \(P_i=\X/G_i\), the (definable) set of orbits of the \(G_i\) action on 
  \(\X\), and let \(\pi_i:\X\ra{}P_i\) be the quotient map. We set 
  \(\ti\Cc=\bigcup\{\Cc,P_i\}\). We let \(\QF\) be the collection of 
  definable subsets of Cartesian products of \(\Q\),\(\X\) and \(\ti{\Cc}\) 
  preserved (set-wise) by \(G\). Note that the maps \(\pi_i\) are included in 
  \(\QF\).
  
  Let \(\widehat{G}=\Aut[\ti\Cc]{\Q,\X}\), and we will show that 
  \(G=\widehat{G}\).  We already know that \(\widehat{G}\) is an \ogroup, and 
  \(G\) is contained in \(\widehat{G}\) (in a way compatible with the 
  action.) By definition, \(\widehat{G}\) fixes the \(G\) orbits on \(\X\).  
  However, the action of \(\widehat{G}\) on \(\X\) is free, so 
  \(G=\widehat{G}\).
\end{proof}

\subsection{The opposite groupoid}\label{ssc:H}
Let \(\Q\) be a definable set, \((\X,f)\)-internal to \(\Cc\). We assume the 
internality datum to be given as in the beginning of~\ref{sss:defgp}. In 
particular, we have a definable map \(\pi:\Cc\ra{}D\) whose fibres are image 
sets of elements of \(\X\) (regarded as maps on \(\Q\)), and another 
definable set \(H\) whose elements can be considered as bijective maps 
between the mentioned fibres:
\begin{equation*}
  H=\{\autC{x}{y}\|x,y\in\X\}
\end{equation*}
  
Our purpose in this section is to describe the structure of \(H\).  It turns 
out that \(H\) is a \emph{definable groupoid}\index{groupoid!definable}, 
acting definably on \(\X\).  In particular, we get a family of groups acting 
on \(\X\), and the action turns out to be free.  These groups do not act by 
automorphisms (and in general do not act on \(\Q\) at all), but any point of 
\(\X\) gives rise to a (non-canonical) isomorphism of any of these groups 
with the automorphisms group \(G\) of the internality structure. In 
particular, there is only one isomorphism class of these groups, which is 
determined by any of them.  Given a collection of definable sets \(\QF\), 
these statements go through for \(G_\QF\) and an \(\w\)-groupoid \(H_\QF\).  
The advantage of considering \(H\) and not \(G\) is that \(H\) belongs to 
\(\Se{\Cc}\).

\subsubsection{Definable sets of types}
Recall that the composition of elements of \(\X\) and \(H\), viewed as 
functions from \(\Q\) to \(\Cc\) and from \(\Cc\) to \(\Cc\), is denoted by 
\(\mu:H\x_D\X\ra\Wbar{\X}\), where \(\Wbar{\X}\) is a set containing \(\X\).  
Consider the definable set defined by \(\mu(h,x)\in\X\).  We view it as a 
family of subsets of \(H\) parametrised by \(\X\), and consider the canonical 
family \(d:\ti{H}\MapsTo{}E\) and the natural map \(t:\X\MapsTo{}E\) obtained 
from it.  Note that \(E\) is in \(\Se{\Cc}\), and, by the proof of 
proposition~\ref{prp:autispair}, the fibres of the map \(t\) are the orbits 
of the action of the definable automorphisms group \(G\) on \(\X\). So \(E\) 
is the set of such orbits, and may be regarded as the definable set of types 
of elements of \(X\) over \(H\), with respect to the formula 
\(\mu(h,x)\in\X\).  The map \(t\) then associates to each element its type.

As with \(\Cc\) and \(D\), we replace \(H\) by \(\ti{H}\), and thus we get an 
action \(\mu:H\x_E{}\X\ra\X\). For \(e\in{}E\), we denote by \(H_e\) and 
\(\X_e\) the fibres over \(e\).

If \(h\in{}H_e\) and \(x\in\X_e\), let \(f=t(\mu(h,x))\). We claim that for 
any other \(y\in\X_e\), \(t(\mu(h,y))=f\) as well. In fact, there is an 
element  \(g\in{}G\) such that \(y=g(x)\), hence, since \(G\) acts by 
automorphisms of the internality structure,
\begin{equation*}
  \mu(h,y)=\mu(h,g(x))=g(\mu(h,x))
\end{equation*}
holds as well, so \(\mu(h,x)\) and \(\mu(h,y)\) are in the same orbit. Thus 
any \(h\in{}H_e\) maps \(\X_e\) bijectively to some \(\X_f\). Let 
\(c:H\ra{}E\) be the map assigning to each \(h\in{}H_e\) the above element 
\(f\). Let \({H_e}^f\) be the set of elements \(h\) such that \(d(h)=e\) and 
\(c(h)=f\).

If \(\QF\) is any collection of definable sets, the construction is 
analogous.  \(E=E_\QF\) is again defined to be the set of orbits of the 
action of \(G_\QF\) on \(\X\), where \(G_\QF\) is the automorphism group 
associated with \(\QF\), as described above. The main difference with the 
case of pure internality is that \(G_\QF\) is an \ogroup, rather than a 
definable group, so \(E_\QF\) is a pro-definable set, and each \({H_e}^f\) is 
an \(\w\)-definable set.

The following is just a restatement of the above construction:
\begin{prop}
  Let \(\QF\) be a set of formulas containing \(f\) and the formula 
  \(x\in\X\).  For any \(a\in\X_e\), \(b\in\X_f\) there is a unique 
  \(h\in{}H_e^f\) with \(\mu(h,a)=b\).  Given any \(g\in{}G_\QF\), 
  \(g(b)=\mu(h,g(a))\).
\end{prop}
\begin{proof}
  \(h\) is simply the element corresponding to \(\autC{a}{b}\). Indeed, we 
  have by definition that \(\mu(h,a)=b\), so in particular \(h\in{}H_e^f\).  
  The fact that \(h\) commutes with \(G_\QF\) follows, as before, from the 
  fact that the action of both is given by composition of functions, on 
  different sides.  This also implies uniqueness, since \(G_\QF\) acts 
  transitively on each \(\X_h\).
\end{proof}

\subsubsection{A description in terms of definable groupoids}
Recall, from~\Cite{groupoid}, that a \emph{definable 
groupoid}\index{groupoid!definable} is a collection of definable sets 
satisfying the axioms of a groupoid\index{groupoid}, namely, a category all 
of whose morphisms are isomorphisms.  Explicitly, we have:
\begin{enumerate}
  \item
    Two definable sets \(\Ob\) (objects) and \(\Mor\) (morphisms.)
    
  \item
    Two definable maps \(dom,cod:\Mor\MapsTo\Ob\), the domain and the 
    co-domain (range) of a morphism, giving rise to a combined map 
    \(\Mor\MapsTo\Ob\x\Ob\). The fibre over the objects \(x,y\in\Ob\) is 
    denoted \(\Mor(x,y)\).
    
  \item
    A definable composition map \(m:\Mor\x_\Ob\Mor\MapsTo\Mor\), where the 
    fibre product is with respect to \(dom\) in the first factor, and with 
    respect to \(cod\) in the second one. This map is over \(\Ob\x\Ob\), with 
    the obvious maps.

  \item
    A definable map \(id:\Ob\MapsTo\Mor\), the identity morphism for any 
    object, over \(\Ob\x\Ob\) (with the diagonal map from \(\Ob\).)
\end{enumerate}

All this data satisfies the usual axioms of a category, as well as the axiom 
stating that every element of \(\Mor\) has an inverse with respect to \(m\) 
(so that all morphisms are isomorphisms.)

Analogously to the case of groups, we may define an 
\(\w\)-groupoid\index{groupoid@\(\w\)-groupoid} to be a pro-definable set 
that has a defining system consisting of definable groupoids and functors.

Furthermore, given a definable groupoid \(\G\), a 
\(\G\)-torsor\index{torsor@\(\G\)-torsor} is a definable surjective family 
\(F\ra\Ob\) and a definable action \(\mu:\Mor\x_\Ob{}F\ra{}F\), such that
\begin{equation*}
\vcenter{\xymatrix{
  \Mor\x_\Ob F \ar[rr]^\mu \ar[dr]_{cod\circ\pi_1} & & F \ar[dl] \\
                                         & \Ob &
                                       }}
\end{equation*}
commutes,

\begin{gather}
  \mu\circ(m\x 1)=\mu\circ(1\x\mu):\Mor\x_\Ob\Mor\x_\Ob F\ra F\\
  \pi_2=\mu\circ(id\x 1):\Ob\x_\Ob F\isom F
\end{gather}

(i.e., the action is compatible with composition and identity maps), and such 
that the map
\begin{equation*}
  \vcenter{\xymatrix{
    \Mor\x_\Ob F \ar[rr]^{\mu\x 1} \ar[dr]_{dom\x cod} & & F\x F \ar[dl] \\
                                           & \Ob\x\Ob &
                                         }}
\end{equation*}
is an isomorphism. When \(\G\) is an \(\w\)-groupoid, a \(\G\)-torsor is a 
system of torsors over the groupoids involved in defining \(\G\).

If \(G\) is an \(\w\)-group acting  on a definable set \(F\) over a map 
\(F\ra{}E\), such that each fibre \(F_e\) is a \(G\)-torsor, then \(G\x{}E\) 
can be viewed as an \(\w\)-groupoid \(\G\) (with both \(dom\) and \(cod\) the 
projections), and \(F\) is then a \(\G\) torsor.

Back to our construction, in this terminology we have an \(\w\)-groupoid with 
morphisms set \(H\) and objects set \(E\), and \(\X\ra{}E\) is a torsor over 
it.  The action of \(G_\QF\) on \(\X\) preserves this map and acts on each 
fibre as a torsor, and so gives rise to another groupoid with \(\X\) as a 
torsor. The two torsor structures commute. However, the groupoid \(H\) is 
connected (i.e., for any \(e,f\in{}E\), \(H_e^f\) is non-empty.)

In particular, for any element \(e\in{}E\), \(H_e^e\) is an \(\w\)-group 
(defined over \(e\)), and any element \(x\in\X_e\) gives rise to an 
isomorphism between the groups \(G\) and \(H_e^e\) (sending the element 
\(g\in{}G\) to the element in \(H\) representing the map \(\autC{g(x)}{x}\).) 
This map is not canonical, though: the map determined by a different element 
\(y\) is obtained from the first one by conjugation with \(\autC{x}{y}\).

The \(H\)-torsor structure of \(\X\) allows us to interpret each element of 
\(H\) as a partial definable function from \(\X\) to itself. Conversely, let 
\(h\) be a partial function from \(\X\) to itself, definable in \(T_\QF\) 
over \(\Se{\Cc}\), and let \(x\in\X_e\) be a point where \(h\) is defined, 
\(y=h(x)\). Then \(h\) coincides on \(\X_e\) with \(\autC{x}{y}\). In fact, 
\(h\) commutes with the action of \(G_\QF\) on \(\X\), and \(\X_e\) is a 
\(G\)-torsor, so \(h\) is defined on the whole \(\X_e\), and is determined by 
its value on any one point, for instance \(x\), where it has the same value 
as \(\autC{x}{y}\). Thus for any \(e,f\in{}H\), \(H_e^f\) is naturally 
interpreted as the set of \(T_\QF\) definable functions from \(\X_e\) to 
\(\X_f\) defined over \(\Se{\Cc}\), and the groupoid structure is just 
composition of functions (which turn out, in particular, to be bijections.)

\subsection{Summary}\label{ssc:summary}
We conclude the section with a summary of the notation and results that 
appeared so far.
\begin{notation}\label{not:aut}\index{internality!summary of notation}
Let \(\Q\) be \((\X,f)\) internal to \(\Cc\), \(\QF\) a collection of 
definable sets.
\begin{itemize}
  \item
    \(D\) is the set of image sets of \(\Q\) under elements of \(\X\).  
    Without loss of generality we assume that the subsets of \(\Cc\) 
    corresponding to distinct elements of \(D\) are disjoint. The family of 
    inverses to the family \(f\) is denoted by \(g:\Cc\x_D\X\ra\Q\).

  \item
    \(F\) is the set of bijections from \(\Q\) to itself, obtained by 
    composing two elements of \(\X\) with the same image. \(G=G_\QF\) is the 
    group of automorphisms of \(\Q\) (and \(\X\)), preserving \(\QF\) and the 
    internality datum. In particular, \(G_\emptyset\) is the group of 
    automorphisms of the internality datum, and 
    \(G_\QF\subseteq{}G_\emptyset\). \(F\) need not be a group in general, 
    but \(G_\emptyset\subseteq{}F\), compatibly with the action on \(\Q\).  
    The action of \(G_\emptyset\) on \(\X\) is by composition of functions 
    from \(\Q\) to \(\Cc\).
    
  \item
    \(E=E_\emptyset\) is the set of orbits of the action of \(G_\emptyset\) 
    on \(\X\). We denote by \(t\) the map from \(\X\) to \(E\). Thus we have 
    a ``forgetful'' map from \(E\) to \(D\). \(H\) is the set of partial maps 
    from \(\Cc\) to itself, obtained by composing two elements of \(\X\).  If 
    \(h\in{}H\) is obtained as \(\autC{x}{y}\), composition with \(h\) maps 
    bijectively the \(G_\emptyset\) orbit of \(x\) to that of \(y\). Thus, by 
    modifying \(H\), we get maps from \(H\) to \(E\), for the domain and 
    image of an element of \(H\). The composition map is denoted by 
    \(\mu:H\x_E\X\ra\X\). The composition of the map from \(H\) to \(E\x{}E\) 
    with the map to \(D\) gives the domain and range of elements of \(H\) as 
    maps on \(\Cc\).

    \(H=H_\QF\) and \(E_\QF\) are the analogous constructions for the group 
    \(G_\QF\). Thus \(E_\QF\) is the (pro-definable) set of orbits of \(G\) 
    on \(\X\), and \(H\) is the opposite groupoid to \(G\).
\end{itemize}
\end{notation}

The main results are summarised in the following theorem.
\begin{thm} \label{thm:abstract}
  \mbox{}
  \begin{enumerate}
    \item
      Let \(\Q\) be \((\X,f)\)-internal to \(\Cc\)
      \begin{enumerate}
        \item\label{thi:abstract}
          For any collection \(\QF\) of definable sets, there is an \ogroup 
          \(G_\QF\), and a definable action of \(G_\QF\) on \(\Q\) and 
          \(\X\), whose points in a model \(M\) are identified via this 
          action with the group of automorphisms of \(\Q,\X\) preserving all 
          \(\QF\) sets and the internality structure (and fixing all other 
          sorts involved in \(\QF\).) This group is the intersection of the 
          definable groups \(G_{\QF_0}\) for finite subsets \(\QF_0\) of 
          \(\QF\).

        \item
          If \(G\) is an \ogroup acting faithfully on \(\Q\), then it is of 
          the form mentioned above (perhaps for different \(\X\) and \(f\))

      \end{enumerate}

    \item
      If \(\Q\) is internal to \(\Cc\), \(\QF\) a collection of definable 
      subsets of \(\Q,\Cc\), \(G=G_\QF\) the group as in~\ref{thi:abstract}, 
      and \(\Cc\) is stably embedded, then \(G\) does not depend on the 
      internality datum.

    \item
      If \(\Q\) and \(\Cc\) are definable sets, \(\Cc\) is stably embedded, 
      and for some saturated model \(M\), the automorphism group \(G\) of 
      \(\Q[M]\) over \(\Cc[M]\) has the same cardinality as \(M\), then 
      \(\Q\) is internal to \(\Cc\) (and the internality structure is 
      automatically preserved by \(G\)), and \(G\) is, therefore, an \ogroup.
  \end{enumerate}
\end{thm}
\begin{proof}
  \mbox{}
  \begin{enumerate}
    \item
      \begin{enumerate}
        \item
          This is proposition~\ref{prp:defaut}.

        \item
          proposition~\ref{prp:defisaut}.
      \end{enumerate}

    \item
      This follows from proposition~\ref{prp:intinEQ}, since in this case the 
      internality datum is constructed from \(\Q\) and \(\Cc\) using \(\QF\).

    \item
      This is proposition~\ref{prp:card}.\qedhere
  \end{enumerate}
\end{proof}


\section{Stable theories with a generic automorphism}\label{sec:stability}

We consider a theory \(T_\sg\)\index{aa@\(T_\sg\)}, whose models are models 
of a given theory \(T=\EQ{T}\) endowed with an automorphism \(\sg\), which is 
generic\index{automorphism!generic}, in a sense defined below. We will 
consider internality datum in \(T_\sg\), where the set \(\Cc\) will be the 
set of fixed points of \(\sg\). Our goal, which we achieve under some 
additional assumptions, will be to describe, in terms of \(T\), the 
automorphism group preserving all the \(T\) structure. More precisely, we 
will describe \(\QF\)-types in \(\X\) in terms of \(T\) definable invariant 
functions on \(\X\).  Our main application is the case where \(T\) is the 
theory \ACF{} of algebraically closed fields, which is dealt with in 
section~\ref{sec:ACFA}.

\subsection{The theory \(T_\sg\)}
Let \(T\) be an arbitrary theory that eliminate quantifiers. Let \(B\) be a 
definably closed subset of a model of \(T\), and let \(\sg_0\) be an 
automorphism of \(B\) (i.e., a bijection of \(B\) with itself, preserving all 
the quantifier free relations.) We denote by \(B_0\) the subset of \(B\) 
consisting of elements fixed by \(\sg_0\). We will consider only models of 
\(T\) that contain \(B_0\), so we assume that \(B_0\) is contained in 
\(\Dcl(0)\) (if \(B_0\) is non-empty, it follows that \(T\) is complete.)

Let \(M\) be any model of \(T\) containing \(B\), \(\sg\) an automorphism of \(M\) 
extending \(\sg_0\), and \(A\subseteq{}M\) a definably closed subset of \(M\), 
closed under \(\sg\) and containing \(B\). We call such a pair \((A,\sg)\) a 
\emph{\(\sg\)-structure}\index{structure@\(\sg\)-structure}. The theory of 
\emph{\(T\) with a generic automorphism} is defined to be a theory 
\(T_\sg\)\index{aa@\(T_\sg\)|(} in the language 
\(L_{\sg}=L\cup{}{\sg}\cup{}B\) (where \(\sg\) is a unary function symbol) 
with the properties that:
\begin{itemize}
  \item
    \(T_\sg\) contains the universal theory of \(T\), and says that \(\sg\) is a 
    map from the universe to itself that is a homomorphism of the \(T\) 
    structure (i.e., preserves all sets definable in \(T\). In particular, it 
    is injective.)

  \item
    Every \(\sg\)-structure can be embedded in a model of \(T_\sg\). This means 
    that for any \(\sg\)-structure \((A,\sg)\) there is a model \(N\) of \(T_\sg\) 
    and a function \(f:A\ra{}N\) over \(B\), such that \(f\circ\sg=\sg_N\circ{}f\), 
    and for any tuple \(\Wbar{a}\in{}A\), any formula over \(B\) satisfied by 
    \(\Wbar{a}\) in a model of \(T\) containing \(A\) is also satisfied by 
    \(f(\Wbar{a})\) in \(N\) (this does not depend on the model of \(T\) used, since 
    \(T\) is model complete.)

  \item
    \(T_\sg\) itself is model complete.
    
\end{itemize}

In other words, \(T_\sg\) is the \emph{model companion} (cf~\Cite{sacks}) of 
the theory of \(T\) with an arbitrary automorphism. \(T_\sg\) need not exist, 
in general, but if it exists, it is unique. The existence and properties of 
such theories were studied in~\Cite{PillayZoe}. It follows from the second 
condition (with \(\Wbar{a}\) the empty tuple) that \(T\) coincides with 
\(T_\sg\) restricted to \(L\), that any model of \(T_\sg\) is a model of 
\(T\) together with an automorphism, and that any definable set of \(T\) can 
be identified with a definable set of \(T_\sg\).\index{aa@\(T_\sg\)|)}

We may now state our assumptions (which will be valid till the end of the 
section):
\begin{Ass}\label{ass:stable}
  We are given a theory \(T\) in a language \(L\), which is assumed to 
  eliminate quantifiers, and to have EI. We fix a structure \(B\), and an 
  automorphism \(\sg_0\) of \(B\) as above. We assume:
  \begin{enumerate}
    \item
      \(T_\sg\), the theory of \(T\) with a generic automorphism as described 
      above, exists. This theory will play the role of the general theory 
      \(T\) in section~\ref{sec:general}.

      We denote by \(T_\sg^0\) the theory obtained from \(T_\sg\) by 
      restricting the constant symbols to \(B_0\).

    \item
      \(T\) is a stable\index{stable} theory

    \item
      \(T_\sg^0\) eliminates imaginaries.
      
    \item \label{asi:stableint}
      \emph{Assumptions about the internality}:

      \(\Cc\) is the set \(\sg(x)=x\) of fixed points of \(\sg\).

      \(\Q\) is a \(T_\sg\) definable set, internal to \(\Cc\). We further 
      assume that this internality is witnessed by a \(T_\sg\) definable set 
      \(\X\) that is given within a \(T\) definable set \(\ti{\X}\) by the 
      formula \(\sg(x)=A(x)\), where \(A\) is a \(T_B\) definable map.

      Moreover, we assume that the maps \(\pi:\X\ra{}D\), 
      \(g:\Cc\x_D\X\ra\Q\) and \(\mu:H\x_E\X\ra\X\) are given by terms 
      (function symbols) in \(T_\sg\).

    \item
      \(\QF\) is the collection of quantifier free sets in \(T_\sg\).
  \end{enumerate}
\end{Ass}

\begin{remark} \label{rmk:qfstar}
  \mbox{}
  \begin{enumerate}
    \item
      The EI assumption for \(T_\sg\) is discussed in~\Cite{groupoid}. It 
      is shown there that this condition can be translated to a condition on 
      \(T\), namely, that in \(T\) there are no nontrivial definable groupoids 
      with finite \(Hom\) sets. There is also a description of a procedure for 
      adding sorts to \(T\) to achieve this condition, similar to the way this 
      is done to obtain EI.

    \item\index{stability}
      All our use of stability is concentrated in two results, 
      claim~\ref{clm:acl} and claim~\ref{clm:cb}. The background from 
      stability theory required for these results is explained in 
      appendix~\ref{ssc:stable}.

    \item\label{rmi:qf}
      It follows from the fact that \(\pi\), \(g\) and \(\mu\) are given by 
      terms (item~\ref{asi:stableint} of assumption~\ref{ass:stable}) that 
      composition with them does not increase the number of quantifiers.  In 
      particular, it follows that the collection \(\QF^*\) mentioned 
      in~\ref{sss:defaut} (and in the proof of proposition~\ref{prp:defaut}) 
      is the set of quantifier free sets. It is this condition that we 
      actually use.
  \end{enumerate}
\end{remark}

We conclude immediately from EI in \(T_\sg^0\) that \(\Cc\) is stably 
embedded:
\begin{lemma}
  Assume~\ref{ass:stable}. Then \(\Cc\) is stably embedded (in \(T_\sg\))
\end{lemma}
\begin{proof}
  Let \(\phi(c,x,b)\) be a family of subsets of \(\Cc\) parametrised by 
  \(x\), with \(b\in{}B\) and \(\phi\) in \(T_\sg^0\). By EI, there is a 
  \(T_\sg^0\) canonical family \(\psi(c,z)\) and function \(h\) such that 
  \(\phi(c,x,y)\iff{}\psi(c,h(x,y))\). These sets are preserved by \(\sg\), 
  and since \(\sg\) fixes \(\Cc\) pointwise, it fixes the canonical 
  parameter, so the image of \(h\) is contained in \(\Cc\).  Now, \(\psi\) 
  and \(h(-,b)\) give a canonical family for the original set.
\end{proof}

\subsubsection{The definability of \(\sg(x)\)}
We make a few comments on the seemingly strong condition that for \(x\in\X\), 
\(\sg(x)\) is definable (in \(T_B\)) over \(x\). We first note that it is 
enough to require that for \emph{some} \(x\in\X\), \(\sg(x)\) is definable 
over \(x\), since we may restrict to the subset given by such a definition.

By taking ``prolongations'', we may replace this condition by the condition 
that \(\sg^n(x)\) is definable over
\begin{equation*}
x,\sg(x),\dots,\sg^{n-1}(x)
\end{equation*}
for some \(n\). This is because \(\X\) can then be replaced by the subset 
\begin{equation*}
(x_0,\dots,x_{n-1})\in\X\x\sg(\X)\x\dots\x\sg^{n-1}(\X)
\end{equation*}
given by \(\sg(x_i)=x_{i+1}\). This set is definably isomorphic to \(\X\), 
and therefore all assumptions are preserved. It satisfies the  requirement 
that \(\sg(\Wbar{x})\) is definable over \(\Wbar{x}\), and since \(\Cc\) is 
stably embedded, the automorphism group does not depend on the internality 
datum.

We next show that this condition does in fact hold under each of two 
assumptions. The first is a further assumption about the shape of the 
internality datum. We recall that in general we have (using the notation as 
in~\ref{ass:stable}) a definable group action \(m:G_\emptyset\x\X\ra\X\) with 
a definable quotient map \(\pi_E:\X\ra{}E\), where \(E\) is contained in 
\(\Cc^n\). We also have an action \(\mu:H\x_E\X\ra\X\) of the groupoid \(H\) 
on \(\X\), and the two actions commute. In general, this structure is defined 
over \(B\). The assumption of the next proposition is essentially that there 
is an extension to a similar structure over \(B_0\). More precisely, we have:
\begin{prop}
  Let \(G_\emptyset\) be the definable automorphism group of the internality 
  datum, let \(E\) be the set of orbits of the action 
  \(m:G_\emptyset\x\X\ra\X\) of \(G_\emptyset\) on \(\X\), and let 
  \(\pi_E:\X\ra{}E\) be the quotient map (so \(x\Conj{}y\) if and only if 
  \(\pi_E(x)=\pi_E(y)\)). Assume that there is a \(T_\sg\) definable action 
  \(m_1:G_1\x\X_1\ra\X_1\), such that
  \begin{enumerate}
    \item \(G_\emptyset\le{}G_1\), \(\X\subseteq\X_1\) and \(m_1\) extends 
      \(m\)
    \item The quotient \(\pi_1:\X_1\ra{}E_1\) of this action extends 
      \(\pi_E\) (i.e., if \(m_1(g,x)=y\) where \(x,y\in\X\), then 
      \(g\in{}G_\emptyset\)).
    \item \(\mu\) extends to a map \(\mu_1:H\x_{E_1}\X_1\ra\X_1\), and 
      \begin{equation*}
        m_1(g,\mu_1(h,x))=\mu_1(h,m_1(g,x))
      \end{equation*}
      for all \(g\in{}G_1\) and \((h,x)\in{}H\x_{E_1}\X_1\).
    \item \(\X_1\), \(\pi_1\) and \(\mu_1\) are defined over \(B_0\)
    \item \(m_1\) is the restriction of a definable set in \(T_B\).
  \end{enumerate}
  Then for \(x\in\X\), \(\sg(x)\) is definable over \(x\).
\end{prop}
\begin{proof}
  We first note that we may assume that \(E_1\subseteq\Cc^k\) for some \(k\).  
  This is because we know this for \(E\), and we may replace \(E_1\) with 
  \(E_1\cap\Cc^k\), and \(X_1\) with its inverse image.

  Now let \(x\in\X\) be any element. Then \(x\in\X_1\), and since \(\X_1\) is 
  over \(B_0\), also \(\sg(x)\in\X_1\). Also, since \(\pi_1\) is over 
  \(B_0\), and since \(E_1\) is in the constants, we get 
  \(\pi_1(\sg(x))=\sg(\pi_1(x))=\pi_1(x)\). Therefore, there is a unique 
  \(A_x\in{}G_1\) such that \(m_1(A_x,x)=\sg(x)\).
  \begin{claim}
    \(A_x\) does not depend on \(x\)
  \end{claim}
  \begin{proof}
    Let \(y\in\X\) be any other element, and let \(h\in{}H\) be the element 
    corresponding to \(\autC{x}{y}\). Then \(\mu_1(h,x)=y\). Applying \(\sg\) 
    to this equation (and using the assumption that \(\mu_1\) is over 
    \(B_0\)), we get \(\sg(y)=\sg(\mu_1(h,x))=\mu_1(h,\sg(x))\). Hence
    \begin{equation*}
      \sg(y)=\mu_1(h,\sg(x))=\mu_1(h,m_1(A_x,x))=m_1(A_x,\mu_1(h,x))=m_1(A_x,y)
    \end{equation*}
    By uniqueness, \(A_x=A_y\).
  \end{proof}
  Let \(A\) be the constant value of all \(A_x\). \(A\) is thus defined over 
  \(0\), and for all \(x\in\X\) we have \(\sg(x)=m_1(A,x)\).
\end{proof}

The second variant is an assumption on \(T\). A sub-structure \(A\) of a 
model \(M\) of \(T\) is finitely generated over another subset \(B\) if 
\(A\subseteq\Dcl(a,B)\) for a finite tuple \(a\in{}A\).
\begin{prop}
  Assume that \(T\) satisfies the condition that for any structure \(B\), 
  every sub-structure of a structure that is finitely generated over \(B\) is  
  itself finitely generated. Then for some \(n\),  \(\sg^n(x)\) is definable 
  (in \(T\)) over \(\{\sg^i(x)\|0\le{}i<n\}\).
\end{prop}
\begin{proof}
  Let \(B_1\) be the \(\sg\)-structure generated by an element \(b\in\X\) (so 
  this is the \(T\)-structure generated by \(\sg^i(b)\), for integer \(i\).)
  If \(a\in\X\) is some other element, let \(h\in{}H\) be the element 
  taking \(b\) to \(a\). Then the \(\sg\)-structure \(A_1\) generated by 
  \(a\) over \(B_1\) is contained in the \(T\)-structure generated by \(h\) 
  over \(B_1\) (since \(\sg\) is the identity on \(h\), and \(B_1\) is closed 
  under \(\sg\).) Hence, by assumption, \(A_1\) is finitely generated over 
  \(B_1\) as a \(T\) structure. In particular, for some \(m\), \(\sg^m(a)\) 
  is \(T\) definable over \(\sg^i(a)\) for \(i<m\), and \(B_0\).  However, 
  \(b\) was an arbitrary element, so it may be taken to be independent (in 
  the sense of stability) from \(a\). It follows (see claim~\ref{clm:free}) 
  that \(\sg^n(a)\) is \(T\) definable from the \(\sg^i(a)\), \(i<m\).  
  Applying \(\sg\) enough times to the definition, we may also get 
  \(i\ge{}0\).
\end{proof}

We note that the assumption on finitely generated structures is true when 
\(T\) is \(\w\)-stable. In fact, let \(D\) be a sub-structure of the 
structure generated by some tuple \(c\). It is enough to show that there is 
no infinite strictly increasing chain of sub-structures of \(D\).  However, 
the function assigning to each sub-structure \(D_0\) the Morley rank and 
degree of \(c\) over \(D_0\) is strictly decreasing, since for 
\(d\in{}D\setminus{}D_0\), \(d\in{}\Dcl(c)\setminus{}D_0\), hence the Morley 
rank and degree of \(d\) over \(c\) is strictly smaller than the same over 
\(D_0\), so the result follows by symmetry.

\subsection{Invariant functions}\label{ssc:invariant}
Given a \(T_B\) definable relation \(h(\Wbar{x})\), we denote by \(h^\sg\) 
the relation \(h(\inv{\sg}(\Wbar{x}))\), so that, for any model \(M\), 
\(h^\sg(M)=\sg(h(M))\) (note that this is again \(T_B\) definable.) In 
particular, if \(h\) is a \(T_B\) definable function on \(\X\) (that is, a 
\(T_B\) definable relation whose restriction to \(\X\) is a function), then 
\(h^\sg\) is the function on \(\sg(\X)\) obtained by conjugation with 
\(\sg\): \(h^{\sg}=\sg\circ{}h\circ\inv\sg\). For such a function \(h\) we 
also denote by \(h_A\) the function on \(\X\) given by composition with the 
\(T_B\) definable function \(A\) above: \(h_A=h^\sg\circ{}A\). We call a 
\(T_B\) definable function \(h\) on \(\X\) \emph{invariant}\index{invariant 
function} if \(h=h_A\) as functions on \(\X\) (so that for \(x\in\X\), 
\(h(x)\in\Cc\).)

Under the above assumptions, we shall prove:

\begin{prop}\label{prp:stable}
  Assuming~\ref{ass:stable}, let \(M\) be any \(\sg\)-structure, 
  \(a,b\in{}\X(M)\). Then the relation
  \begin{equation*}
    \tp[\Cc(M)]{a}=\tp[\Cc(M)]{b}
  \end{equation*}
  is given by the set of formulas \(h(a)=h(b)\), where \(h\) is a \(T_B\) 
  definable  invariant function. In particular, it is \(\QF\)-definable 
  (i.e., defined by an intersection of quantifier free formulas in \(T\), 
  possibly infinite in number).
\end{prop}

Together with the explicit description in~\ref{sss:defaut}, this implies the 
following description of the automorphism group, which is the main result of 
this section:

\begin{cor} \label{cor:stable}
  Assume~\ref{ass:stable}. The group \(G=\Aut[\Cc]{\Q}\) is given by the 
  intersection of formulas of the form
  \begin{align*}
    G_h(g)=\forall x\in\X(h(x)=h(gx))
  \end{align*}
  where \(h\) is a \(T_B\) definable invariant function.
\end{cor}

\begin{proof}
  By proposition~\ref{prp:defaut} and the remarks following it, the group is 
  given by the intersection of formulas saying that \(x\) and \(gx\) have the 
  same \(\QF^*\) type over \(\Cc(M)\cup{}H(M)\cup{}D(M)\), where \(x\in{}M\).  
  We first note that \(D\), \(H\) and \(E\) are canonical families of subsets 
  of \(\Cc\), and therefore are subsets of \(\Cc\) themselves, since \(\Cc\) 
  is stably embedded.

  Thus we are reduced to equality of types over \(\Cc(M)\). As mentioned in 
  remark~\ref{rmk:qfstar}.\ref{rmi:qf}, \(\QF\) is preserved by composition 
  with the functions \(g\) and \(\mu\) appearing in the description of 
  \(\QF^*\).  Hence the corollary follows from proposition~\ref{prp:stable}.
\end{proof}

We use the terms \emph{algebraic closure} and \emph{definable closure} to 
mean these concepts with respect to \(T\), whereas we say \emph{\(\sg\)-algebraic 
closure}, \emph{\(\sg\)-definable closure} for the same concepts in \(T_\sg\). 
Similarly we write \(\Acl\), \(\Dcl\) in \(T\) and \(\Acl_\sg\), \(\Dcl_\sg\) in 
\(T_\sg\).

The proof of proposition~\ref{prp:stable} depends on two general claims, 
given below. Both involve types in infinitely many variables. These are 
simply maximal consistent sets of formulas in these variables (over a given 
set.) We call such types \emph{infinitary} for short. Any subset \(B\) of a 
model has a type over any other set \(A\), just like in the finite case, which 
we denote as usual by \(\Tp{B}{A}\) (the variables of this type will be indexed 
by the elements of \(B\); thus a statement of the form \(\Tp{B}{A}=\Tp{C}{A}\) 
implies we are given a bijection between \(B\) and \(C\).) A realisation is also 
defined in the same way as for types of finite tuples.

\begin{claim} \label{clm:acl}
  Let \(T\) be a stable theory, let \(T_\sg\) be its associated theory with a 
  generic automorphism, and let \(A\) be an algebraically closed 
  \(\sg\)-structure. Then \(A\) is \(\sg\)-algebraically closed.
\end{claim}

\begin{claim} \label{clm:cb}
  Let \(M\) be a model of a stable theory \(T\) with EI, \(C\subseteq{}M\) a 
  definably closed subset, \(A,B\subseteq{}M\). Let \(E_A\) be the set of 
  elements of \(M\) fixed by all automorphisms that fix \(A\) pointwise, and 
  fix \(C\) as a set, \(C_A=E_A\cap{}C\) (and similarly for \(B\).)
  If \(\Tp{A}{C_A\cup{}C_B}=\Tp{B}{C_A\cup{}C_B}\). Then 
  \(\Tp{A}{C}=\Tp{B}{C}\)\footnote{
    Note that if \(C\) is the set of \(M\) points of a \emph{definable} set, 
    then it is automatically preserved by any automorphism, and \(E_a\) (in a 
    saturated model) is simply \(dcl(a)\). In this case, the condition thus 
    says that \(C\) is stably embedded, so the claim implies that any 
    definable set in a stable theory is stably embedded. If \(C\) is not 
    definable, this equivalence can be viewed as the \emph{definition} of 
    being stably embedded (and this claim says that \emph{any} definably 
    closed subset in a model of a stable theory is stably embedded)}
\end{claim}

Claim~\ref{clm:acl} is a generalisation of the same result for \ACFA{}, as 
appears in~\Cite{ACFA}. It was also proven in~\Cite{PillayZoe}. The second 
essentially appears in~\Cite[Ch.~7,remark~1.16]{Pillay}. Both of these claims 
are consequences of stability, and are explained in 
appendix~\ref{ssc:stable}.  Meanwhile, we use them to deduce 
proposition~\ref{prp:stable}:

\begin{Proof}[of proposition~\ref{prp:stable}]
  Let \(\ti{M}\) be a model of \(T_\sg\), extending \(M\). The collection of 
  \(T\)-definable functions into \(\Cc\) will remain the same in \(\ti{M}\).  
  Also, since \(M\) is definably closed, the  values of any \(T_B\) definable 
  invariant function \(h\) on \(a\) and \(b\) lies in \(\Cc(M)\). Therefore, 
  we may assume that \(M\) is a model of \(T_\sg\). We write \(\Tp{A}{B}\) 
  for the type in the sense of \(T\).

  Let \(a_1\) and \(a_2\) be elements of \(\X(M)\), such that for any \(T_B\) 
  definable invariant function \(h\), \(h(a_1)=h(a_2)\). We would like to 
  show that \(\Tp{a_i}{\Cc(M)}\) are equal in \(T_B\). However, since \(B\) 
  may contain elements not fixed by \(\sg\), it is more convenient to 
  consider these elements as additional parameters. Therefore, we set 
  \(B_i=\Dcl(B\cup{}a_i)\). The map that is the identity on \(B\) and takes 
  \(a_1\) to \(a_2\) extends canonically to a bijection from \(B_1\) to 
  \(B_2\), and in the sense of  this map, we shall prove that 
  \(\Tp{B_1}{\Cc(M)}=\Tp{B_2}{\Cc(M)}\).  This will be done by proving the 
  same for an increasing sequence of subsets of \(\Cc(M)\), until we arrive 
  at the assumptions of \ref{clm:cb}. We note that since \(\sg(a_i)=A(a_i)\), 
  and \(B\) is preserved by \(\sg\), so is each \(B_i\).
  
  We first note that the assumptions on \(a_i\), together with elimination of 
  imaginaries imply that \(\Tp{B_1}{0}=\Tp{B_2}{0}\) (recall that \(0\), the  
  definable closure of the empty set in \(T\), does not include any non-fixed 
  parameters): The canonical parameter of any definable  set lies in \(\Cc\), 
  and the function taking an element of the set to this parameter is 
  invariant and \(B\) definable, hence \(a_1\) and \(a_2\) agree on this set.

  By definition, \(B_i=\Dcl(B\cup{}a_i)\) is the image of all \(B\)-definable 
  functions on \(a_i\). If \(h\) is a \(T_B\) definable function and 
  \(\sg(h(a_i))=h(a_i)\), we may pass to a globally invariant function \(t\)  
  on \(\X\) (as in the beginning of~\ref{ssc:invariant})
  by defining \(t(x)\) to be \(h(x)\) if \(h(x)=h_A(x)\), and a constant 
  fixed element outside this set. We will have \(h(a_i)=t(a_i)\).
  Hence \(\Cc(B_i)\) is the image of all invariant \(T_B\) definable 
  functions on \(a_i\).
  
  Therefore, the assumption on \(a_i\) implies that \(\Cc(B_1)=\Cc(B_2)\). We 
  denote this set by \(D\).  Further, the same assumption means that 
  \(\Tp{B_1}{D}=\Tp{B_2}{D}\): If \(\phi(a_1,b,d)\) holds for some 
  \(b\in{}B\), \(d\in{}D\), then \(d=h(a_1)=h(a_2)\) for some invariant 
  \(h\), and so \(\phi(a_1,b,h(a_1))\) is a \(B\)-definable formula that 
  holds of \(a_1,b\), hence also of \((a_2,b)\) by the previous step. We next 
  show that the same holds when \(\Cc(B_i)\) is replaced by 
  \(\Cc(\Acl[B_i])\).

  Indeed, for any \(c\in\Cc\), \(\Tp{c}{B_i}\) (hence its set of 
  realisations) is preserved by \(\sg\), since \(\sg(c)=c\) and \(B_i\) is a 
  \(\sg\)-structure. In particular, if \(c\) is algebraic over \(B_1\), the 
  finite set \(Y\) of conjugates of \(c\) over \(B_1\) is fixed by \(\sg\).  
  Since \(T\) has EI, this set is coded by an element \(y\in\Cc\), definable 
  over \(B_1\). By assumption, it is definable over \(B_2\), and \(B_1\), 
  \(B_2\) have the same type over \(y\).  Hence \(c\) is algebraic over 
  \(B_2\), and \(B_1\), \(B_2\) have the same type over \(c\): if \(c\) 
  satisfies some formula over \(B_1\), so do all its conjugates over \(B_1\), 
  hence all its conjugates over \(B_2\). Since \(E=\Cc(\Acl[B_1])\) is the 
  set of such \(c\), the same is true for this set.

  As we just showed, \(E=\Cc(\Acl[B_2])\) and \(\Tp{B_1}{E}=\Tp{B_2}{E}\).  
  We now use again the fact that each \(B_i\) is a \(\sg\)-structure, and 
  deduce from claim~\ref{clm:acl} that \(\Acl\) can be replaced by 
  \(\Acl_\sg\), so that \(E=\Cc(\Acl_\sg(B_i))\).
  
  Finally, we note that the condition \(\tp{a_1}=\tp{a_2}\) does not depend 
  on the model \(M\). Therefore, we may assume that \(M\) is saturated.  
  Under this assumption, the set \(E_{B_i}\) that appears in 
  claim~\ref{clm:cb} is contained in \(\Dcl_\sg(B_i)\), and so the set 
  \(C_{B_i}\) there is contained in \(E\). Since we just proved that 
  \(\Tp{B_1}{E}=\Tp{B_2}{E}\), and \(C_{B_1}\cup{}C_{B_2}\subseteq{}E\), 
  applying claim~\ref{clm:cb} we get the result.
\end{Proof}
 
Remark~\ref{rmk:cb} ties the description of the images of invariant functions 
with canonical bases.

\begin{remark}
  As can be seen from the proof, internality is not used in this proposition, 
  but only the assumption that \(\sg(a)\) is definable (in \(T_B\)) over 
  \(a\).
\end{remark}

\subsubsection{Comparison with the classical group}
We now apply these results again, in order to compare the group \(G\) of 
automorphisms of the quantifier free structure, to the usual model theoretic 
automorphism group.  Let \(G_0=Aut(\Q/\Cc)\) be the subgroup of \(G\) 
preserving all \(\Cc\)-definable subsets of \(\Q\). It turns out that \(G_0\) 
is pretty close to \(G\):

\begin{prop} \label{prp:profinite}
  Assume~\ref{ass:stable}, and let \(G=\Aut[\Cc]{\Q}\) be the quantifier free 
  automorphism group, and \(G_0\) be the subgroup of full automorphisms.  
  Then the quotient \(G/G_0\) is pro-finite.
\end{prop}

\begin{proof}
  Let \(Y\) be a \(\Cc\) definable subset of \(\X\) (possibly with 
  quantifiers.) We will show that its orbit (as a set) under \(G\) is finite.  
  This will be enough, since \(G_0\) is the stabiliser, inside \(G\), of such  
  subsets.

  Let \(Z\subseteq{}\X\) be an orbit of \(G\). For any point \(z\in{}Z\), we 
  get from the internality datum a subset 
  \(Y_z=\{f_z\circ{}\inv{f_y}|y\in{}Y\}\) of \(H\). Since \(\Cc\) is stably 
  embedded, the canonical parameter \(c_z\) for this set lies in \(\Cc\).  
  Thus we get a \(\sg\)-definable function from \(Z\) to \(\Cc\), sending 
  \(z\) to \(c_z\). By claim~\ref{clm:acl}, \(c_z\) is algebraic (in the 
  sense of \(T\)) over the \(\sg\)-structure generated by \(z\), and 
  therefore over \(z\) (since \(\sg(z)\) is \(T\)-definable over \(z\)). In 
  other words, we have a \(T\) definable function \(t\) from \(Z\) to finite 
  subsets of \(\Cc\), such that \(c_z\in{}t(z)\).  But since \(Z\) is an 
  orbit of \(G\), \(t\) is constant on \(Z\).  Thus there is a finite subset 
  \(W\), such that \(c_z\in{}W\) for all \(z\), so the number of sets \(Y_z\) 
  is finite.  Finally, we note that by fixing an element \(z_0\) of \(Z\), 
  the sets \(Y_z\) are identified with the orbit of \(Y\) under \(G\): any 
  \(g\in{}G\) can be written as \(\inv{f_{z_0}}\circ{}f_{y_g}\) for a unique 
  \(y_g\in{}Z\), hence \(gY=\inv{f_{z_0}}Y_{y_g}\).
\end{proof}

\section{The case of \ACFA{}}\label{sec:ACFA}
In this section we study the example of the situation in~\ref{ass:stable} 
where \(T\) is the theory of algebraically closed fields. The resulting 
theory \(T_\sg\) is called the theory of algebraically closed fields with an 
automorphism (\ACFA{}), and was studied in~\Cite{ACFA}. In particular, it is 
proved there that \ACFA{} eliminates imaginaries\index{elimination of 
imaginaries!for \ACFA}. The base set \(B\) is, in this case, a field with an 
automorphism, which we shall also denote by \(\kk\). An interesting example 
is when \(\kk=\Q(t)\), and \(\sg(t)=t+1\).

\subsection{Internality in \(ACFA\)}
We proceed to describe explicitly the internality datum and the associated 
structure in the example we consider (linear difference equations over 
\(\kk\)). \index{internality!datum!for a linear difference equation} We use 
the notation as in~\ref{ssc:summary}.

The set \(\Cc\), given by the equation \(\sg(x)=x\) is a pseudo-finite field, 
called the fixed field. The set \(\Q\), in our example, is given by an 
equation \(\sg(q)=Aq\), where \(q\) is a (column) tuple of variables (of 
length \(n\)), and \(A\) is an invertible matrix over \(B\). Such an equation 
is called a linear difference equation. The set \(\Q\) has a definable vector 
space structure over \(\Cc\), of dimension \(n\).  Therefore, it is internal 
to \(\Cc\), with the internality datum consisting of the set of vector space 
bases \(\X\), where \(f:\Q\x\X\ra\Cc^n\) assigns to any vector \(q\) and 
basis \(x\) the coefficients of the representation of \(q\) in the said 
basis. We think of the elements of \(\X\) as matrices, whose columns are the 
basis elements (hence solution to the equation.) In these terms, \(f\) is 
given by \(f(q,x)=\inv{x}q\).

Thus, the image of any element \(x\in\X\) (viewed as a map from \(\Q\) to 
\(\Cc^n\) via \(f\)) is the whole \(\Cc^n\), so the set \(D\) consists of one 
point. The inverse map \(g:\Cc^n\x\X\ra\Q\) is given by \(g(a,x)=xa\).  The 
set \(\X\) coincides with the subset of \(\Gl{n}\) given by \(\sg(x)=Ax\).  
The set \(H\) is the set \(\{\inv{y}x\|x,y\in\X\}\), or equivalently, the 
subset of \(\Gl{n}\) given by \(\sg(x)=x\), and so, in other words, is 
identified with \(\Gl{n}(\Cc)\). The action \(\mu:H\x\X\ra\X\) is given by 
\(\mu(h,x)=x\inv{h}\), so the set \(E=E_\emptyset\) also consists of one 
point. Finally, the group \(G_\emptyset\) of automorphisms of the internality 
datum coincides with the set \(F\), defined in section~\ref{ssc:Internality} 
as \(\{y\inv{x}\|x,y\in\X\}\). It is the subgroup of \(\Gl{n}\) given by 
\(\sg(x)=Ax\inv{A}\).

It follows from this description that all assumptions of~\ref{ass:stable} are 
satisfied. The fact that \(\QF\) is the set of quantifier free sets means 
that we are interested in the group \(G\) of automorphisms of \(\Q\) 
preserving all \emph{polynomial} relations. Since \ACFA{} does not have 
quantifier elimination, this is different, in general, from the usual model 
theoretic group \(G_0\) that preserves all definable relations:

\begin{example}\label{xpl:profinite}
  Assume \(B=\QQ\), and let \(\Q\) be given by \(\sg(x)=4x\). It is easy to 
  see that any non-zero solution to this equation is transcendental over 
  \(\Cc\) (in general, in dimension \(1\) the only equations over the fixed 
  field that can have algebraic solutions are those with \(A\) a root of 
  unity. This can be seen by considering the minimal polynomial of a 
  solution, and also follows from proposition~\ref{prp:GgenA}, below.) 
  Therefore, the quantifier free automorphism group is the multiplicative 
  group of \(\Cc\). In particular, there is only one quantifier free type of 
  non-zero elements of \(\Q\) over \(\Cc\). On the other hand, a square root 
  of such a non-zero element may satisfy either the equation \(\sg(x)=2x\) or 
  \(\sg(x)=-2x\) (and the other one satisfies the same equation), so this 
  unique quantifier free type splits into at least two full types.
\end{example}

Note, however, that according to proposition~\ref{prp:profinite}, \(G_0\) is a 
pro-finite index subgroup of \(G\).

We next note that \(\X\) is ``Zariski dense'' in the set \(\ti{\X}\) of bases of 
\(K^n\), i.e., if \(p\) is any polynomial on \(\ti\X\) (over any base set) such 
that \(\sg(x)=Ax\implies{}p(x)=0\), then \(p\) is identically \(0\). This follows 
from the axioms of \ACFA{}. The main axiom of \ACFA{} states 
(cf.~\Cite{ACFA}):

\begin{axiom}\label{axm:ACFA}
  Let \(U\) be an irreducible variety (over the model \(K\)), and let 
  \(V\subseteq{}U\x\sg(U)\) be an irreducible sub-variety, projecting 
  dominantly to each factor. Then for any proper closed subset \(W\) of \(V\), 
  there is a point \(x\in{}U(K)\) with \((x,\sg(x))\in{}V\setminus{}W\).
\end{axiom}

Applying this axiom with \(U=\sg(U)=\ti\X\), \(V\) the variety given by \(Y=AX\) 
and \(W\) the closed subset given by \(p(X)=0\) (and using the fact that \(\ti\X\) 
itself is dense in the \(n^2\) dimensional affine space), we get the result.

It follows that the generic type over \(B\) (in the sense of \ACF{}) is 
consistent with \(\X\). Hence, for the definition of the group we may 
concentrate on types that extend this generic type. This allows us to 
describe the automorphism group \(G\):

\begin{cor}\label{cor:acfa}
  The\index{automorphisms group!first order definition of!for linear 
  difference equations} group \(G\) of automorphisms of the equation 
  \(\sg(q)=Aq\) preserving all polynomial identities, is the subgroup of 
  \(\Gl{n,B}\) given by the equation \(\sg(g)=Ag\inv{A}\), and all the equations
  \begin{equation*}
    \forall x(r(x)\neq 0 \implies (r(gx)\neq 0\land h(x)=h(gx)))
  \end{equation*}
  with \(h=\frac{p}{r}\) a rational function on \(\ti\X\) over \(B\), satisfying 
  (globally) \(h=h^\sg\circ{}A\), where \(h^\sg\) is the function obtained from 
  \(h\) by applying \(\sg\) to the coefficients.

  In particular, it is the intersection of the definable group given by 
  \(\sg(g)=Ag\inv{A}\) and an algebraic group \(\ti{G}\) defined over \(B\).
\end{cor}
\begin{proof}
  Any such rational function \(h\) restricts to an actual function around any 
  point satisfying \(r(x)\neq{}0\). This function has values in \(\Cc\) for 
  elements in \(\X\). If \(g\in{}G\) is any automorphism, then 
  \(r(x)\neq{}0\) implies \(r(gx)\neq{}0\), and therefore \(h(x)=h(gx)\) for 
  any element in \(x\in\X\) with \(r(x)\neq{}0\). Since the condition is a 
  polynomial equation, it is satisfied for any other \(x\) as well. Thus any 
  automorphism satisfies these equations.

  Conversely, let \(g\) be any element satisfying the above formulas. To show 
  that it is an automorphism, it is enough to show for one element \(x\in\X\) 
  that \(x\) and \(gx\) satisfy the same quantifier free type over \(\Cc\).  
  Let \(x\) be an element of \(\X\) generic (in the sense of \ACF{}) over 
  \(g\) (i.e., \(x\) does not lie in any \(g\) definable sub-variety.) Then 
  both \(x\) and \(gx\) are generic over \(0\), and in particular 
  \(r(x)\neq{}0\) and \(r(gx)\neq{}0\).  Any invariant definable function (as 
  in corollary~\ref{cor:stable}) coincides on the generic type with an 
  invariant rational function. Hence, \(x\) and \(gx\) agree on all invariant 
  functions, so by proposition~\ref{prp:stable}, they have the same type.

  The last statement follows since any \(\w\)-definable group in an 
  \(\w\)-stable theory (like \ACF{}) is, in fact, definable 
  (cf.~\Cite{stability}), and any definable group in \ACF{} is algebraic 
  (cf.~\Cite{HrWeil}.)
\end{proof}

\begin{remark}\label{rmk:types}
  We thus have the following description of the set of types: there is an 
  open subset \(U\) of \(\ti\X\) (given by the functions \(r\) above), and an 
  algebraic map \(F\) into some affine space \(L\). The space of types of 
  elements of \(\X\cap{}U\) is the image of \(\X\cap{}U\) under this map, and 
  the image lies in \(L(\Cc)\).

  On the other hand, since \(\ti\X\) is a torsor over \(\Gl[B]{n}\), and 
  \(\ti{G}\) is an algebraic subgroup of \(\Gl[B]{n}\), the quotient 
  \(\ti\X/\ti{G}\) is an algebraic variety \(V\) (over \(B\).) Since the 
  group \(\ti{G}\) is, in general, not preserved by \(\sg\), \(\sg\) is not 
  defined on this quotient. However, if we define \(\phi\) as 
  \(\inv{A}\circ\sg\) on \(\ti\X\) and by \(\phi(g)=\inv{A}\sg(g)A\) on 
  \(\Gl{n}\), we an automorphism of the group action of \(\Gl{n}\) on 
  \(\ti\X\) that preserves the group \(\ti{G}\). Therefore, it does induce an 
  automorphism of \(V\). The sets \(\X\) and \(G\) are precisely the sets of 
  fixed points of \(\phi\), and therefore the set of types \(\X/G\) embeds 
  via the quotient map into the set of fixed points of \(\phi\) on \(V\).

  Hence, if we take the open set \(U\) above to be closed under the 
  \(\ti{G}\) action, we get a map from the image of \(U\) in \(V\) to \(L\), 
  which is a bijection on the level of the sets of types. However, we do not 
  know whether this map can be taken to be an isomorphism, and whether it can 
  be extended to the whole quotient space \(V\).
\end{remark}

In the case that the base field \(B\) consists only of fixed elements (i.e., 
\(\Cc[B]=B\)), we may give a more explicit description of the situation. In 
this case \(\sg\) is, itself, an automorphism, and therefore the 
corresponding element \(A\) of \(\Gl{n}\) belongs to \(G\). In general, we 
claim:

\begin{prop}\label{prp:GgenA}
  Let \(G_A\) be the intersection of all algebraic subgroups of \(\Gl{n}\) 
  defined over \(\Cc(B)\) and containing \(A\), and let \(U\) be the definable 
  subgroup of \(\Gl{n}\) given by \(\sg(x)=Ax\inv{A}\). Then 
  \(G\subseteq{}G_A\cap{}U\). If \(B=\Cc(B)\), the groups are equal.
\end{prop}

As will be evident from~\ref{ssc:singer}, this is similar to 
\Singer{Proposition~1.21}.

\begin{proof}
  Let \(G_0=G_A\cap{}U\). Since \(G_0\subseteq{}U\), \(G_0\) acts on \(\X\).  
  Let \(S_0=\X/G_0\) and \(S=\ti\X/G_A\). Since \(G_A\) is defined over the 
  fixed field, \(\sg\) is well defined on \(S\). Since \(U\) is precisely the 
  subgroup of \(\Gl{n}\) preserving \(\X\), \(S_0\) embeds into \(S\).  The 
  image of \(S_0\) in \(S\) is fixed pointwise by \(\sg\), since 
  \(A\in{}G_A\). Hence elements with the same (quantifier free) type over 
  \(\Cc\) are in the same \(G_A\) orbit. Since these types are \(G\)-torsors, 
  \(G\subseteq{}G_0\).

  If \(B=\Cc[B]\), then \(A\in{}G\), so in particular it belongs to the 
  algebraic group \(\ti{G}\) associated with \(G\). Since this group is also 
  defined over \(B=\Cc[B]\), we get that \(G_A\subseteq{}\ti{G}\).
\end{proof}

\begin{remark}
  If \(E\) is a connected algebraic subgroup of \(\Gl{n}\) containing \(A\), 
  we may also find a solution to the equation in \(E\) (when \(\ti\X\) and 
  \(\Gl{n}\) are identified as algebraic varieties.) This again follows from 
  axiom~\ref{axm:ACFA}: The projections from the subset of \(E\x{}E\) given 
  by \(y=Ax\) to each of the components is an isomorphism.
\end{remark}

\subsubsection{The opposite groupoid}
The situation with the opposite group (studied in section~\ref{ssc:H}) is 
slightly more complicated. According to proposition~\ref{prp:stable}, there 
is a stratification of \(\ti\X\), such that the situation described in 
remark~\ref{rmk:types} holds, maybe with different data on each stratum. To 
simplify notation, we will now assume that the whole of \(\ti\X\) is one 
stratum. We thus have an algebraic map \(F\), with \(E_\QF=F(\X)\).

Given a value \(e\) in this set \(E_\QF\), we obtain a (quantifier free) type 
\(p_e\) of an element of \(\X(M)\), over \(\Cc(M)\) (which is the orbit 
associated with this value), and a subgroup \(H_e\), whose action (given by 
\((x,g)\mapsto{}x\inv{g}\)) is transitive on the realisations of \(p_e\) 
(recall that, in the language of groupoids, \(e\) is an object, \(H_e\) is 
the group of automorphisms of \(e\), and \(p_e\) is the set of isomorphisms 
between \(e\) and the special object corresponding to the equation. The 
action is given by composition).

\begin{prop}
  The group \(H_e\) is given explicitly by the formulas \(\sg(g)=g\) and the 
  algebraic subgroup \(\ti{H_e}\) of \(\Gl{n}\) defined by
  \begin{equation}\label{eqn:He}
    \forall{}x(F(x)=e{}\iff{}F(x\inv{g})=e)
  \end{equation}
\end{prop}
\begin{proof}
  By definition, \(H_e\) is the subgroup of \(\Gl{n}(\Cc)\) given by
  \begin{equation*}
    \forall{}x\in\X(F(x)=e{}\iff{}F(x\inv{g})=e)
  \end{equation*}
  Hence, the only thing that should be verified is that 
  ``\(\forall{}x\in\X\)'' can be replaced by ``\(\forall{}x\)''.
  
  Let \(Z\) be the irreducible algebraic set given by the equation 
  \(F(x)=e\), viewed over the subfield \(L\) of \(M\) generated by \(\kk\) 
  and \(\Cc(M)\). Then \(\X\cap{}Z\) is dense in \(Z\), since, by the 
  definition of \(F\) (\(e\) is the canonical base of \(p_e\)), the \(ACF\) 
  type of any element of \(\X\cap{}Z\) over \(L\) is the generic type of 
  \(Z\). Hence if \((X\cap{}Z)g=X\cap{}Z\) then \(Zg=Z\).
\end{proof}

In general, there is no reason for these subgroups \(H_e\) to coincide (see 
example~\ref{xpl:distinctH}.) However, as explained in section~\ref{ssc:H}, 
they are all conjugate over \(\Cc\), together with the torsors \(p_e\) on 
which they act. We would like to obtain some information about the conjugacy 
classes. Since \(\ti{H_e}\) is the Zariski closure of \(H_e\) (over \(\Cc\)), 
we get that the \(\tilde{H_e}\) are conjugate as well.

We have three notions of conjugacy between the elements of the family 
\(\ti{H_e}\): conjugacy over \(\Cc\), conjugacy over \(\Acl(\Cc)\) 
(``absolute''), and conjugacy over \(\Cc(\kk)\). We will discuss only the 
first two. Given any definable family of subgroups, conjugacy of two elements 
over \(\Cc\) is a definable property, so it makes sense to ask in which field 
a conjugacy class lies (the absolute conjugacy class lies in some 
pro-definable set).

The conjugacy class over \(\Cc\) belongs to \(\Dcl_\sg(\kk)\), since it is 
canonically associated with the equation, and to \(\Cc\), since it is defined 
in terms of subgroups of \(\Cc\) (we now consider the group alone, without the 
torsor.) The absolute conjugacy class is coarser, and therefore belongs to 
the same set. We now make the following observation:
\begin{claim}\label{clm:disjoint}
  For any difference sub-field \(\kk\) of a model \(M\) of \ACFA{}, \(\kk\) is 
  linearly disjoint from \(\Cc(M)\) over \(\Cc(\kk)\). In particular,
  \begin{equation*}
    \Dcl_\sg(\kk)\cap\Cc(M)\subseteq{}\Acl(\Cc(\kk))\cap\Cc(M)
  \end{equation*}
\end{claim}
\begin{proof}
  We need to show that any subset of \(\Cc(M)\) linearly dependent over 
  \(\kk\) is linearly dependent over \(\Cc(\kk)\). Let \(\sum{}c_ia_i=0\) be 
  a minimal linear dependence, where \(c_i\in\Cc(M)\) and \(a_i\in\kk\).  
  Applying \(\sg\), we get \(\sum{}c_i\sg(a_i)=0\). The minimality implies 
  that there is \(a\in\kk\) such that \(\sg(a_i)=aa_i\) for all \(i\). It 
  follows that \(\sg(\frac{a_i}{a_j})=\frac{a_i}{a_j}\) for all \(i,j\).  
  Dividing the dependence by \(a_1\) we thus get a linear dependence over 
  \(\Cc(\kk)\).

  The ``in particular'' part is now deduced as follows: We saw that
  \begin{equation*}
    \Dcl_\sg(\kk)\subseteq\Acl_\sg(\kk)=\Acl(\kk)
  \end{equation*}
  Hence it is enough to prove that
  \begin{equation*}
    \Cc(\Acl(\kk))=\Cc(\Acl(\Cc(\kk)))
  \end{equation*}
  This follows from the fact that \(\Cc(\Acl(\kk))\) is linearly disjoint 
  from \(\kk\) over \(\Cc(\kk)\).
\end{proof}

In particular, if the fixed field of the base is algebraically closed, this 
fixed field will contain a point in the family \(\ti{E_\QF}\) (the Zariski 
closure of \(E_\QF\)), which will therefore define a group over the base 
field that lies in the same (absolute) conjugacy class as any \(H_e\) (note 
that in this case, by the above claim, the conjugacy class will belong to the 
base field.)

\subsection{The algebraic theory}\label{ssc:singer}
Our aim now is to describe the relation between our results and the algebraic 
Galois theory of linear difference equations. This theory is described 
in~\Cite{singer} (and more recently, in~\Cite{SingerMod}.) For the sake of 
clarity, we shall repeat part of the exposition there.

As before, we have a fixed base field \(\kk\), with a fixed automorphism \(\sg\) 
on it. A \emph{difference algebra} over \(\kk\) is a \(\kk\)-algebra \(B\), 
together with a ring automorphism of \(B\) whose restriction to \(\kk\) is 
\(\inv{\sg}\). The inverse of \(\varphi\) will be denoted \(\sg\). A map of 
difference algebras is a map of usual \(\kk\) algebras that commutes with 
\(\sg\). The set of difference algebra maps from \(A\) to \(B\) will be denoted by 
\(\DHom{A}{B}\). A difference ideal is an ideal \(I\) of \(B\) such that 
\(\sg(I)\subseteq{}I\). If \(B\) is Noetherian (which it is, in any situation we 
consider), this implies that \(\sg(I)=I\).  The kernel of a map of difference 
algebras is a difference ideal, and the quotient of a difference algebra by a 
difference ideal is again a difference algebra.

Given a linear difference equation \(\sg(\Wbar{x})=A\Wbar{x}\), where, as before, 
\(\Wbar{x}\) is a tuple of variables, and \(A\) is a matrix over \(\kk\) (which will 
be fixed from now on), we associate with it a difference algebra \(R\) as 
follows: as a ring, \(R=\kk[X,\inv{\Det{X}}]\), where \(X\) is a matrix of 
variables (of the same size as \(A\)), and \(\Det{X}\) is the determinant 
polynomial in these variables (in other words, \(R\) is the coordinate ring 
over \(\kk\) of the variety of invertible matrices.) The action of \(\sg\) on \(R\) 
is determined by the action on the generators, where it acts as \(\sg(X)=AX\).

If \(B\) is any difference algebra (always over \(\kk\)), the set \(\DHom{R}{B}\) 
may be identified with the set of invertible matrices \(x\) over \(B\) that 
satisfy \(\sg(x)=Ax\). In particular, if \(M\) is a model of \(ACFA_\kk\), we have 
\(\X(M)=\DHom{R}{M}\). For any element \(x\in\X(M)\), we will denote by 
\(\phi_x:R\ra{}M\) the corresponding map. The kernel of such a map is a 
difference ideal. This ideal clearly depends only on the quantifier free type 
of \(x\) over \(\kk\), and in turn determines this type (together with the 
difference equation.)

More generally, Let \(L\) be any subfield of \(\Cc(M)\) extending 
\(\Cc(\kk)\), and let \(R_L=L\tensor_{\Cc(\kk)}R\) (note that 
\(\kk_L=L\tensor_{\Cc(\kk)}\kk\) is again a field, since, by 
claim~\ref{clm:disjoint}, \(\kk\) and \(L\) are linearly disjoint over 
\(\Cc(\kk)\); it is the subfield of \(M\) generated by \(\kk\) and \(L\).  
Therefore \(R_L\) is the analogue of \(R\) for the base field \(\kk_L\).) As 
before, any element \(x\in\X(M)\) determines a map \(\phi_x^L:R_L\ra{}M\) and 
an ideal \(I_x^L\) in \(R_L\) that corresponds to the (quantifier free) type 
of \(x\) over \(L\).

Recall also that with any such element \(x\) we have associated the 
\Emph{canonical base} of its type over \(\Cc(M)\), which is a certain 
subfield \(L_x\) of \(\Cc(M)\) (since it is definably closed.) This is, in 
fact, the field generated by the values on \(p\) of the invariant functions.  
The connection between these objects is given by the following proposition.

\begin{prop}
  Let \(p\) the quantifier free type of a solution of the equation over 
  \(\Cc(M)\), let \(I\) be the ideal corresponding to \(p\) in 
  \(R_{\Cc(M)}\), and let \(D=R_{\Cc(M)}/I\).
  \begin{enumerate}
    \item
      Let \(L\) be a perfect subfield of \(\Cc(M)\), \(I^L=I\cap{}R_L\) the 
      ideal corresponding to \(p\) in \(R_L\), \(D_L=R_L/I^L\). Then the map 
      \(\Cc(M)\tensor_{R_L}D_L\ra{}D\) is an isomorphism if and only if \(L\) 
      contains the canonical base of \(p\).

    \item
      The ideal \(I\) is a maximal difference ideal in \(R_{\Cc(M)}\).
  \end{enumerate}

  In particular, if \(L\) contains the canonical base of \(p\), then \(I^L\) 
  is a maximal difference ideal.
\end{prop}
\begin{proof}
  \mbox{}
  \begin{enumerate}
    \item
      If \(L\) contains the canonical base of \(p\), \(I^L\) contains all 
      elements of the form \(f(x)-c\) where \(f\) is an invariant function.  
      By proposition~\ref{prp:stable}, this determines the type \(p\) 
      completely.  Hence \(I\) is generated by \(I^L\).

      Conversely, if the map is an isomorphism, let \(d\) be the value on 
      \(p\) of an invariant function \(f\). Then \(d\tensor{}1-1\tensor{}f\) 
      goes to \(0\) in \(D\), so it is \(0\) already in 
      \(\Cc(M)\tensor_{R_L}D_L\). Hence \(d\in{}L\). Since the canonical base 
      of \(p\) is the definable closure of all these values, and \(L\) is 
      definably closed, \(L\) contains the canonical base.
    \item
      Assume that \(I\) is not maximal, and let \(J\supset{}I\) be a 
      difference ideal extending \(I\). Since \(I\) comes from a solution 
      that lies in the field \(M\), it is a prime ideal. Therefore, the 
      dimension of \(J\) is strictly smaller than the dimension of \(I\). We 
      may also assume that \(J\) is radical, since the radical of a 
      difference ideal is also a difference ideal. We show that in fact, we 
      may assume that \(J\) is prime. Let \(J_1,\dots,J_n\) be the prime 
      decomposition of the ideal generated by \(J\) in \(R_{\Acl(\Cc(M))}\), 
      as a usual ideal (so that \(J=\bigcap{}J_i\).) This decomposition 
      happens, in fact, over some finite field extension of \(\Cc(M)\), and 
      in particular, over the fixed field of some \(\sg^m\). Then \(\sg^m\) 
      acts on this set of ideals, and for some larger \(m\), \(\sg^m\) fixes 
      each of them.

      We now replace \((M,\sg)\) by \((M,\sg_1)\), where \(\sg_1=\sg^m\). In 
      this new model, the fixed field \(C_1\) in the new structure is an 
      algebraic (finite) extension of the original fixed field. The solution 
      \(x\) of the original equation is also a solution of the implied 
      equation for \(\sg_1\). The associated ideal \(I_1\) extending \(I\) is 
      again a prime difference ideal (for the difference algebra obtained 
      from the new equation.) So is each of the ideals \(J_i\) extending 
      \(I_1\). Since the base field extension is algebraic, the dimensions of 
      \(I\) and \(I_1\) are the same, and similarly for \(J_i\). In 
      particular, each \(J_i\) is a proper extension of \(I_1\). Replacing 
      \(J\) with \(J_i\), this shows that we may assume that \(J\) is prime.

      Since \(J\) is a prime difference ideal, \(R_{\Cc(M)}/J\) is a 
      \(\sg\)-structure (i.e., a difference algebra which is an integral 
      domain.) Therefore, it can be embedded (over \(\Cc(M)\)) in \(M\). If 
      \(q\) is the quantifier free type of the image \(X\) under this 
      embedding, the associated ideal is \(J\). Now let \(x\) be a solution 
      of \(p\), \(y\) a solution of \(q\), \(h=\inv{x}y\). Then \(h\) defines 
      an automorphism of \(R_{\Cc(M)}\) taking \(I\) to \(J\). The image of 
      \(J\) under this map will be strictly contained in \(J\).  This process 
      gives an infinite increasing chain of ideals, contradicting the fact 
      that \(R_{\Cc(M)}\) is Noetherian.\qedhere
  \end{enumerate}
\end{proof}

Note that if \(L\) satisfies the above condition, then \(R_L/{I_x^L}\) is a 
simple difference ring whose fixed field \(L\) is the fixed field of 
\(\kk_L\).  This is similar to \Singer{Lemma~1.8}.

Let \(p\) be a fixed type over \(\Cc(M)\), \(L\) its canonical base, \(I\) 
the associated ideal in \(R_L\) and \(S_L=R_L/I\). Thus, the set \(p(M)\) of 
realisations of \(p\) in \(M\) is identified with \(\DHom[\kk_L]{S_L}{M}\).  
Therefore,

\begin{equation}
\DHom[\kk_L]{S_L\tensor_{\kk_L}R_L}{M}=p(M)\x\X(M)
\end{equation}

Recall from section~\ref{ssc:Internality}, that given an element 
\(x\in{}p(M)\), and another element \(y\in\X(M)\), we may form the element 
\(h=\inv{x}y\), which conversely determines, together with \(x\), the element 
\(y\). In other words, any element \(x\in\X(M)\) gives a definable bijection 
between \(\X\) and \(H\) (which in our case is just \(\Gl{n}(\Cc)\).) This is 
reflected by the fact (\Singer{equation~1.2}) that \(S_L\tensor_{\kk_L}R_L\) 
is isomorphic to \(S\tensor_L{}T\), where \(T=L[Y,\inv{\Det{Y}}]\) represents 
\(H\) (i.e., \(H(M)=\DHom[L]{T}{M}\).)

Finally, dividing by the ideal generated by \(I\) in 
\(S_L\tensor_{\kk_L}R_L\), we get the ring \(S_L\tensor_{\kk_L}S_L\), that 
represents the set \(p(M)\times{}p(M)\) of pairs of realisations of \(p\) in 
\(M\). The definable bijections mentioned above restrict to definable 
bijections between \(p\) and the group \(H_p\) we obtained. Since the group 
\(H_p\) is algebraic and defined over \(\kk_L\), it is represented by some 
quotient \(W_L\) of \(T\). The bijection mentioned above corresponds to the 
equation \(S_L\tensor_{\kk_L}S_L=S_L\tensor_L{}W_L\), \Singer{equation~1.3}.  
When \(\Cc(k)\) is algebraically closed, these algebras are, as explained 
earlier, isomorphic (as usual algebras) to algebras \(S\tensor_{\Cc(\kk)}L\) 
and \(W\tensor_{\Cc(\kk)}L\), where \(S\) and \(W\) are defined over \(\kk\).  
This recovers (up to isomorphism) the algebras constructed in~\Cite{singer}.

Note that in contrast with the algebraic approach, there need not be a torsor 
of solutions defined over the base field \(\kk\). This corresponds to the 
cases when the definable set \(E_\QF\) of all torsors (the image of \(\X\) 
under the meromorphic invariant function \(F\)), contains no point of 
\(\kk\). This can happen even if \(\Cc(\kk)\) is algebraically closed, since 
\(E_\QF\) is not (in general) constructible. In this case, we still obtain a 
group and a torsor over it by taking a point in some constructible set 
containing \(E_\QF\). This torsor is isomorphic (over an extension field) to 
the torsors we obtained, and is isomorphic to the torsors obtained 
in~\Cite{singer}. However, both are isomorphisms of algebraic varieties, and 
the points of this torsor do not solve the equation.

\subsection{Some examples}\label{ssc:examples}
\begin{example}
  Consider the equation \(\sg(x)=-x\). Then \(\X\) is the same set, with 
  \(0\) removed. For any \(x\in\X\), \(x^2\) lies in \(\Cc\), and so \(x^2\) 
  is an invariant function. If \(\kk\) contains a solution \(d\) of the 
  equation, then \(x/d\) is also an invariant function, and so the group is 
  trivial.
  
  We now assume that this is not the case. Then we may take \(F(x)=x^2\), and 
  \(E_\QF=F(\X)\) is then the set of fixed elements whose square root is not 
  fixed.  Any of the torsors we defined is of the form \(x^2=e\) with 
  \(e\in{}E_\QF\). The canonical base for this torsor is the field extension 
  of \(\Cc(\kk)\) obtained by adding \(e\). If \(\Cc(\kk)\) has no extensions 
  of order \(2\) (in particular, if \(\Cc(\kk)\) is algebraically closed), 
  then \(E_\QF\) has no \(\kk\) points, and so none of these torsors is 
  defined over \(\kk\).  The torsors considered in~\Cite{singer} are of the 
  form \(x^2=d\) where is a non-zero element of \(\Cc(\kk)\). If \(d\) is not 
  in \(E_\QF\), this set is isomorphic algebraically (over \(\Cc(\kk)(e)\)) 
  to the set \(x^2=e\), but the \(\sg\) structure is not the same.
\end{example}

\begin{example}\label{xpl:distinctH}
  Let \(A\) be of the form
  \begin{equation*}
  \begin{pmatrix}
    -1 & a \\
     0 & b
  \end{pmatrix}
  \end{equation*}
  We may take a matrix of solutions to lie in the set \(\Mat{x}{y}{0}{z}\), 
  where \(x\) satisfies \(\sg(x)=-x\). In other words, the set of triangular 
  matrices is consistent with the set of solutions. However, if we would like 
  to emulate the construction of the automorphisms group \(G\) using global 
  invariant functions (as above), restricting to this set would be 
  counter-productive: the set of upper triangular matrices is a proper closed 
  subgroup, and indeed, the function \(x^2\) is an invariant function on this 
  set, but not on the whole set of solutions.
  
  Instead, let \(Y=\Mat{x}{y}{z}{w}\) be an arbitrary solution.  Apply the 
  determinant, we see that
  \begin{gather*}
    \sg(\Det{Y})=\Det{A}\Det{Y}=-b\Det{Y} \\
    \intertext{Therefore,}
    \sg(\frac{w}{\Det{Y}})=\frac{\sg{w}}{\sg(\Det{Y})}=-\frac{w}{\Det{Y}}
    \intertext{and similarly}
    \sg(\frac{z}{\Det{Y}})=-\frac{z}{\Det{Y}}
  \end{gather*}

  Hence the function \(F\) given by
  \begin{equation*}
    Y\mapsto 
    (\frac{z^2}{\Det{Y}^2},\frac{zw}{\Det{Y}^2},\frac{w^2}{\Det{Y}^2})
  \end{equation*}
  is an invariant function from \(\ti\X\) onto the subset \(\ti{E}\) of 
  \(\A^3\setminus{}(0,0,0)\) given by \(s^2=rt\), and from the whole set of 
  solutions to \(\ti{E}(\Cc)\).
  
  According to proposition~\ref{prp:GgenA}, the automorphism group \(G\) will 
  be contained in the group of matrices of the form 
  \(\Mat{\pm{}1}{x}{0}{y}\), and for generic \(a\) and \(b\) will be equal to 
  this group. For certain values of \(a\) and \(b\), there will be some other 
  invariant functions. For example, if \(\kk=\Cc(\kk)\), \(G\) will be the 
  subgroup generated by \(A\) (again, by proposition~\ref{prp:GgenA}.) The 
  function \(((b+1)x-az)^2\) will be an additional invariant function in this 
  case.
  
  If \(Y_1\) and \(Y_2\) are two solutions with \(F(Y_1)=F(Y_2)\), the matrix 
  taking \(Y_1\) to \(Y_2\) is in the above group. Therefore, for generic 
  \(a\) and \(b\) these three functions generate all other invariant 
  functions. Note that the automorphism group can be obtained, in this case, 
  by imposing the condition that an element preserves any one of the above 
  functions alone, but on the whole space of solutions. However, this will 
  not suffice when restricting to a particular type. Note also, regarding 
  remark~\ref{rmk:types}, that \(\ti{E}\) is the quotient of the action of 
  the group of matrices of the form \(\Mat{\pm{}1}{x}{0}{y}\) on \(\ti{X}\).

  For similar reasons, the groups in the family \(H\) will be distinct in 
  this case. Let \(T_{(d,e,f)}\) be the set of solutions 
  \(Y=\Mat{x}{y}{z}{w}\) corresponding to the point 
  \((d,e,f)\in{}\ti{E}(\Cc)\). Let \(H_{(d,e,f)}\) be the group of matrices 
  that preserve \(T_{(d,e,f)}\) (this is a typical element of the family of 
  groups \(H\).) Then a direct calculation shows that \(H_{(d,e,f)}\) is 
  given by the equations
  \begin{gather*}
    d\frac{p^2}{\Det{Z}^2}+2e\frac{pr}{\Det{Z}^2}+f\frac{r^2}{\Det{Z}^2}=d\\
    d\frac{pq}{\Det{Z}^2}+2e\frac{ps+rq}{\Det{Z}^2}+f\frac{rs}{\Det{Z}^2}=e\\
    d\frac{q^2}{\Det{Z}^2}+2e\frac{qs}{\Det{Z}^2}+f\frac{s^2}{\Det{Z}^2}=f
  \end{gather*}
  where \(Z=\Mat{p}{q}{r}{s}\) is an element of \(H_{(d,e,f)}\). In 
  particular, it does depend on the point \((d,e,f)\).
\end{example}

\begin{example}
  The equality of groups in proposition~\ref{prp:GgenA} is false in general 
  (when \(A\) is not over the fixed field.) For example, consider the 
  one-dimensional equation given by \(A=\frac{t+1}{t}\), where \(\kk=\Q(t)\), 
  with \(\sg(t)=t+1\).  Since \(A\) is not a root of unity, it generates the 
  whole of \(G_m\). However, the automorphism group is trivial, since \(t\) 
  is a solution.
\end{example}


\appendix
\section{Model theoretic background}\label{sec:app}
\subsection{elimination of imaginaries}\label{ssc:EI}
The notions of imaginaries\index{imaginaries}, and elimination of 
imaginaries\index{elimination of imaginaries|(}, were introduced 
in~\Cite{EI}.

Recall that a theory \(T\) has the property of \emph{elimination of 
imaginaries (EI)}\index{EI|see{elimination of imaginaries}} if for any 
definable family \(\phi(x,y)\), there is a definable set \(Z\) and a 
definable family \(\psi(z,y)\) with parameter variable \(z\) in \(Z\), such 
that
\begin{equation*}
  \forall{}x\in{}X\exists!z\in{}Z\forall{}y(\phi(x,y)\iff\psi(z,y))
\end{equation*}
holds in \(T\).
By definition, such \(Z\) and \(\psi\) determine a unique definable map 
\(f_\phi:X\MapsTo{}Z\), such that \(\phi(x,y)\iff{}\psi(f_\phi(x),y)\). If we 
require that \(f_\phi\) is onto (i.e., any member of the \(\psi\) family is 
also a member of the \(\phi\) family), the triple \((Z,\psi,f_\phi)\) is 
determined up to a unique definable map, and is called a \Emph{canonical 
family} for \(\phi\).

If \(\phi(x,y)\subseteq{}X\x{}X\) happens to be a definable equivalence 
relation on \(X\), then \(f_\phi:X\onto{}Z\) as above is simply the quotient 
of \(X\) by this relation. Conversely, given any 
\(\phi(x,y)\subseteq{}X\x{}Y\), the definable set \(E_\phi(x_1,x_2)\) given 
by \(\forall{}y(\phi(x_1,y)\iff\phi(x_2,y))\) is an equivalence relation on 
\(X\), and a canonical family for \(\phi\) as above is the quotient of \(X\) 
by \(E_\phi\). Thus, a theory has EI if and only if any equivalence relation 
has a quotient.

We say that a collection of definable sets \(\{X_i\}\) has EI if the theory 
of these sets with the induced structure (i.e., the theory whose sorts are 
the \(X_i\) and whose basic relations are all the definable subsets of 
products of the \(X_i\) in the original theory) eliminates imaginaries.

This notion depends on the language. In a context where the language can be 
changed, EI can always be assumed. To show this, we modify the language 
as follows:
\begin{itemize}
  \item For any family \(\phi(x,y)\subseteq{}X\x{}Y\) parametrised by \(X\), 
    we introduce a new sort \(X_\phi\) and a new function symbol 
    \(\pi_\phi:X\ra{}X_\phi\).
    
  \item For any such family \(\phi\), let \(\psi(z,y)\subseteq{}X_\phi\x{}Y\) 
    be given by the formula
    \begin{equation*}
      \exists{}x\in{}X(\phi(x,y)\land{}\pi_\phi(x)=z)
    \end{equation*}

    We extend the theory by the requirement that \((X_\phi,\psi,\pi_\phi)\) is 
    a canonical family for \(\phi\).
\end{itemize}

The theory obtained in this way (when applying the procedure for all sorts) 
eliminates imaginaries (it might seem that the process needs to be iterated; 
however, any family in the new theory corresponds to a family in the original 
one, and is therefore accounted for.) This theory is denoted 
\(\EQ{T}\)\index{aa@\(\EQ{T}\)}. The procedure does not change the category 
of models (and elementary maps), up to natural equivalence.  In particular, a 
type over parameters in \(T\) can also be viewed as a type in \(\EQ{T}\). 
However, for a given language and theory, it is an interesting question 
whether the theory in that language admits EI.

Given a collection of definable sets \(X_i\), we denote by 
\(\Eq{X_i}\)\index{aa@\(\Eq{X_i}\)} the collection of all canonical families 
for all families with parameter variables in Cartesian products of the 
\(X_i\) (possibly by adding sorts, as described above.)
A definable map from a definable set \(Y\) to \(\Eq{X_i}\) is a definable map 
into (a finite union of) any of these families (in other words, it is a 
definable set in \(\Eq{X_i,Y}\) that happens to be a function on \(Y\).)
\index{elimination of imaginaries|)}

\subsection{stable embeddedness}\label{ssc:SE}
The notion and properties of stable embeddedness\index{stable embeddedness|(} 
are discussed in the appendix of~\Cite{ACFA}. A definable set \(X\) is called 
\Def{stably embedded} if for any definable family of subsets of \(X^n\) there 
is a family definable with parameter variable in \(X\) and with the same 
fibres. If \(X\) also has EI, this means that for any definable family of 
subsets of \(X\) (with parameter not necessarily in \(X\)), there is a 
canonical family\index{canonical family} with parameter variable in \(X\).  
For example, if \(X\) is an algebraically closed field with no other induced 
structure (but in a theory containing other sets), this means that any family 
of distinct constructible sets is itself constructible.

When \(X\) has EI, the assumption that \(X\) is stably embedded can also be 
stated as follows: for \(a\in{}M\), where \(M\) is some model, the subset 
\(\phi(a,x)\) of \(X(M)\) is determined by the value \(f_\phi(a)\in{}X(M)\) 
of some definable function \(f_\phi\) at \(a\). This description implies that 
\(\Tp{a}{X(M)}\) is determined by its restriction to the values on \(a\) of 
all such functions \(f_\phi\) for all \(\phi\). This is a small set, 
contained in \(\Dcl(a)\). In particular, it does not depend on \(M\) (this 
description does not depend essentially on the EI assumption, since, as 
mentioned above, we may pass to the type over \(\Eq{X}\) in a unique way.)  
For this reason, stably embedded sets enjoy some of the good properties of 
small sets.
\index{stable embeddedness|)}

We conclude these two parts of the appendix with a remark about these notions 
applied to incomplete theories. Though both imaginaries and stable 
embeddedness are defined syntactically, the standard approach in model theory 
is to work in a model. However, once a model is chosen, one is dealing with 
two theories, the original (incomplete) theory, and the theory of the model.  
In the following proposition we the interpretation of several notion in terms 
of \(T\), with the corresponding interpretation in an extension.
\begin{prop}
  Let \(T_1\) be a theory extending \(T\) (in the same language), and let 
  \(M\) be a model of \(T_1\). Any definable set \(X\) of \(T\) determines a 
  definable set in \(T_1\), also denoted by \(X\).
  \begin{enumerate}
    \item If \(T\) eliminates imaginaries, then so does \(T_1\).
    \item If \(X\) is stably embedded in \(T\), then it is stably embedded in 
      \(T_1\)
    \item Assume that \(T\) eliminates imaginaries. If \(A\subseteq{}M\), 
      then any element of \(\Dcl_{T_1}(A)\) is inter-definable in \(T_1\) 
      with an element of \(\Dcl_T(A)\)
    \item The notion of an automorphism of \(M\) does not depend the theory
  \end{enumerate}
\end{prop}
\begin{proof}
  \begin{enumerate}
    \item If \(\phi(x,y)\subseteq{}X\x{}X\) defines an equivalence relation 
      in \(T_1\), the formula \(E_\phi(x,y)\) given by 
      \(\forall{}z(\phi(x,z)\iff\phi(y,z))\) (``\(\phi\) has the same fibre 
      over \(x\) and \(y\)'') gives a definable equivalence relation in 
      \(T\), whose quotient is the same as the quotient of \(X\) by \(\phi\) 
      in \(T_1\).
    \item is obvious
    \item Let \(b\in\Dcl_{T_1}(A)\). Then \(f(a)=b\), with \(a\in{}X(A)\) for 
      some \(T_1\) definable set \(X\) and function \(f:X\ra{}Y\). These can 
      be represented by some definable sets \(\ti{X}\) and \(\ti{Y}\) and a 
      definable relation \(\phi(x,y)\subseteq{}\ti{X}\x\ti{Y}\) in \(T\). If 
      \((\ti{Z},\psi(z,y),\pi)\) is a canonical family for \(\phi\), then 
      \(\psi\) determines a bijective function in \(T_1\) (and hence in 
      \(M\)). Thus, \(b\) is inter-definable with \(\pi(a)\), which is in 
      \(\Dcl_T(a)\).
    \item is also obvious.
  \end{enumerate}
\end{proof}

In particular, the last two parts imply that we may use the usual methods of 
automorphisms to (essentially) determine the definable closure.

\subsection{Stability}\label{ssc:stable}
In\index{stability|(}  this appendix we give a short overview of the notions 
from stability used in the paper. In particular, we give the proofs of 
claim~\ref{clm:acl} and claim~\ref{clm:cb}. The facts below appear in many 
texts on stability, for example~\Cite{stability}.

Let \(T\) be a complete and model complete theory with EI.  Till the end of 
this section the word \emph{set} means a definably closed subset of some 
model of \(T\). For any set \(A\) and any \(A\)-definable set \(X\), we 
denote by \(\Df[A]{X}\) the boolean algebra of \(A\)-definable subsets of 
\(X\), and by \(\Stn[A]{X}\) the space of types over \(A\) that belong to 
\(X\) (if \(A=\Dcl[\emptyset]\), we omit it.) \(T\) is called 
\(\Def{stable}\) if for any algebraically closed set \(A\), any type 
\(p(x)\in\Stn[A]{X}\) is \(A\)-definable.  This\index{definable type} means 
that for any definable set \(Y\), there is a map 
\(d_p:\Df{X\x{}Y}\ra{}\Df[A]{Y}\) such that for any set \(B\supseteq{}A\), 
the set of formulas
\begin{align*}
  \Par{p}{B}\df{}p\cup{}\{\phi(x,b)\|b\in{}B,d_p\phi(b)\}
\end{align*}
is consistent. Note that in this case, by completeness and model completeness 
of \(T\), \(\Par{p}{B}\) is a complete type over \(B\). The collection of 
maps \(d_p\) is called a \Def{definition scheme}. It is a fact that if \(T\) 
is stable, such a definition scheme \(d_p\) is unique.

We now assume that \(T\) is stable. The definition scheme extends uniquely to 
a definition scheme \(d_p:\Df[A]{X\x{}Y}\ra\Df[A]{Y}\) (by changing \(Y\).) 
Each such map \(d_p\) is a homomorphism of boolean algebras.  Therefore it 
induces a map \(d_p^*:\Stn[A]{Y}\ra\Stn[A]{X\x{}Y}\) of type spaces (given by 
the inverse image: \(d_p^*(q)=\{\phi\|d_p(\phi)\in{}q\}\).) So given types 
\(p\in\Stn[A]{X}\), \(q\in\Stn[A]{Y}\) we get a new type 
\(d_p^*(q)\in\Stn[A]{X\x{}Y}\). This type extends both \(p\) and \(q\): For 
\(\phi\in\Df[A]{X}\), \(d_p\phi\) is a sentence over \(A\), which is true 
(and therefore belongs to \(q\)) if and only if \(\phi\in{}p\). On the other 
hand, if \(\phi\in\Df[A]{Y}\), \(d_p\phi=\phi\), so if \(\phi\in{}q\) then 
\(\phi\in{}d_p^*(q)\).

In fact, \(d_p^*(q)\) is the ``freest'' extension of \(p\) and \(q\). This 
can be made precise by (at least) the following claim:

\begin{claim}\label{clm:free}
  Let \(p(x)\), \(q(y)\) be types over an algebraically closed set \(A\) in a 
  stable theory \(T\). If the formula \(\phi(x,y)\df{}(f(x)=g(y))\) belongs 
  to \(d_p^*(q)\), where \(f\), \(g\) are \(A\)-definable functions, then for 
  some \(a\in{}A\), \(f(x)=a\) belongs to \(p\) and \(g(x)=a\) belongs to 
  \(q\).
\end{claim}

\begin{proof}
  The formula \(d_p\phi(y_1)\land{}d_p\phi(y_2)\land{}g(y_1)\ne{}g(y_2)\) is 
  inconsistent: otherwise, for \(y_1\) and \(y_2\) satisfying it \(d_p\) 
  gives an extension of \(p\) containing the formulas \(f(x)=g(y_i)\) and 
  \(g(y_1)\ne{}g(y_2)\), a contradiction.
  
  Therefore, \(g\) is constant on \(d_p\phi\), and since both \(g\) and 
  \(d_p\phi\) are \(A\)-definable, the constant value \(a\) belongs to \(A\).  
  Since \(d_p\phi\in{}q\), this shows that \(g(y)=a\) is in \(q\).  Since 
  \(d_p^*(q)\) extends \(q\), this implies that \(f(x)=a\) is in \(d_p^*(q)\) 
  and therefore in \(p\).
\end{proof}

Given an \(A\)-definable map \(f\), we denote by \(f_*\) (or sometimes by 
\(f\) again) the induced map on types: 
\(f_*(p)=\{\phi(x)\|\phi(f(x))\in{}p\}\). We now note that definition schemes 
are compatible with definable maps: If \(f\) is a definable function on a 
type \(p\) and \(g\) is a definable function on a type \(q\) (all over 
\(A\)), then \((f,g)_*(d_p^*(q))=d_{f_*(p)}^*(g_*(q))\): both are the set of 
formulas \(\phi(x,y)\) such that \(\phi(f(x),g(y))\in{}d_p^*(q)\).  Applying 
this observation to projections, we see that the construction of the type 
\(d_p^*(q)\) extends to types with infinite number of variables: Given two 
such types \(p\) and \(q\), a formula \(\phi(x,y)\) is in \(d_p^*(q)\) if it 
is in \(d_{p_x}^*(q_y)\), where \(p_x\) and \(q_y\) are the restriction \(p\) 
and \(q\) to the variables in \(\phi\).

If \(\sg\) is an automorphism of \(A\) and \(p\) is a (possibly infinite) 
type, we denote by \(\sg(p)\) the type obtained by applying \(\sg\) to all 
the coefficients (so that an extension of \(\sg\) to an elementary map takes 
a realisation of \(p\) to a realisation of \(\sg(p)\).) For a formula 
\(\phi(x,a)\) with \(a\in{}A\), we denote by \(\phi^\sg\) the formula 
\(\phi(x,\inv{\sg}(a))\). 

Now, Given a type \(p(x)\) over \(A\), the map taking a formula \(\phi(x,y)\) 
(over \(\emptyset\)) to \(d'_p(\phi)=(d_{\sg(p)}\phi)^\sg\) is a definition 
scheme for \(p\): Let \(M\) be a model containing \(A\) to which \(\sg\) 
extends. \(b\in{}M\) satisfies \(d'_p(\phi)\) if and only if \(\sg(b)\) 
satisfies \(d_{\sg(p)}(\phi)\).  Therefore, 
\(\sg(p)\cup\{\phi(x,\sg(b))\|d'_p(\phi)(b)\}\) is consistent. Applying 
\(\inv{\sg}\), this shows that \(p\cup\{\phi(x,b)\|d'_p(\phi)(b)\}\) is 
consistent, and so \(d'\) is a definition scheme. By uniqueness, \(d'=d\). 
From all this follows that for any two (possibly infinite) types \(p\), 
\(q\), \(\sg(d_p^*(q))=d_{\sg(p)}^*(\sg(q))\).

Using these notions, we may prove claim~\ref{clm:acl}:
\begin{Proof}[of claim~\ref{clm:acl}]
  Let \(M_1\), \(M_2\) be models of \(T_\sg\) containing \(A\). We will 
  construct a third model \(M\) containing \(M_1\) and \(M_2\) freely over 
  \(A\). We shall then apply this construction in the case \(M_1=M_2\). The 
  freeness will mean that any algebraic element of \(M_1\) outside of \(A\) 
  will give rise to two copies in \(M\). However, the number of elements 
  conjugate to a given algebraic element is known in advance, so this will 
  show that there are no such elements.
  
  Let \(p\) and \(q\) be the types of \(M_1\) and \(M_2\) in \(T\) (over 
  \(A\)), and let \(F=B_1\cup{}B_2\) be a realisation of \(d_p^*(q)\) (we 
  label the variables of this type by the elements of \(F\).) The original 
  automorphisms \(\sg_1,\sg_2\) of the \(M_i\) give rise to a bijective map 
  \(\sg\) from \(F\) to itself. We claim that \(\sg\) is an elementary map.  
  To show this we need to show that the assignment of \(\sg(b)\) to a 
  variable \(x_b\) of \(\sg(d_p^*(q))\) satisfies this type.  However, we 
  just showed that \(\sg(d_p^*(q))=d_{\sg(p)}^*(\sg(q))\), and \(\sg(p)\) is 
  the same type as \(p\) with the variable \(x_{\sg(b)}\) renamed to \(x_b\) 
  (and similarly for \(q\).) Therefore, \(\sg(d_p^*(q))\) is also 
  \(d_p^*(q)\) with the same renaming of variables, so the assignment 
  satisfies this type by the definition of \(F\).
  
  Therefore \(F\) is a \(\sg\)-structure, and by the definition of \(T_\sg\) 
  there is a model \(M\) of \(T_\sg\) containing \(F\). The definition of 
  \(F\) gives rise to an embedding of \(M_1\) and \(M_2\) in \(M\), and since 
  \(T_\sg\) is model complete, this embedding is elementary. Therefore we 
  proved:

  Any two models over \(A\), \(M_1\) and \(M_2\), of \(T_\sg\) can be 
  elementarily embedded over \(A\) into a third model \(M\), such that for 
  any tuples \(x\in{}M_1\) and \(y\in{}M_2\), with \(T_A\) types \(p\) and 
  \(q\), the \(T_A\) type of the pair \((x,y)\) in \(M\) is \(d_{p}^*(q)\).

  Now let \(M\) be any model of \(T_\sg\) and \(a_0\in{}M\) an element in 
  \(\Acl_\sg(A)\). We denote by \(a\) the tuple of conjugates of \(a_0\).  
  Let \(N\) be a model as above, for \(M_1=M_2=M\). Thus \(M\) has two 
  elementary embeddings into \(N\), \(f_1\) and \(f_2\), with \(b_i=f_i(a)\) 
  such that \(\Tp{b_1,b_2}{A}=d_{\Tp{b_1}{A}}^*(\Tp{b_2}{A})\). However, 
  since the embedding is elementary, \(b_1\) and \(b_2\) are both the set of 
  solutions of some \(\sg\)-algebraic formula. Therefore, \(\Tp{b_1,b_2}{A}\) 
  contains the formulas \(\pi_i(b_1)=\pi_j(b_2)\) (where \(\pi_i\) are 
  projections).  By claim~\ref{clm:free}, for any \(i\), 
  \(\pi_i(b_1)\in{}A\), so \(a\in{}A\) as well.
\end{Proof}

We now aim to prove claim~\ref{clm:cb}. We first note that the condition 
\(\Tp{A}{C}=\Tp{B}{C}\) consists of a collection of condition on types of 
finite sub-tuples. Furthermore, for any tuple \(a\), \(E_a\subseteq{}E_A\) 
(in the notation of claim~\ref{clm:cb}.) Therefore, the whole statements 
reduces to the case when \(A\) and \(B\) are finite tuples \(a\) and \(b\).

We explain the result in the case when \(C\) is algebraically closed. In this 
case, let \(p(x)=\Tp{a}{C}\), \(q(x)=\Tp{b}{C}\) and \(C_a\), \(C_b\) as in 
the claim. We first claim that for any \(0\)-definable set \(\phi(x,y)\), 
\(d_p\phi\) is defined over \(C_a\). Indeed, let \(\tau\) be an automorphism 
fixing \(a\) and \(C\) (as a set.) Then \(\tau(p)=p\) and so 
\({(d_p\phi)}^\tau=d_{\inv{\tau}(p)}\phi=d_p\phi\). This shows that 
\(d_p\phi\) is over \(E_a\) and it is over \(C\) by definition.

Now let \(r=\Res{p}{D}=\Res{q}{D}\), where \(D=\Dcl[C_a\cup{}C_b]\). The 
extension of \(r\) to \(\Acl[D]\) is the same in \(p\) and in \(q\): if 
\(d\in\Acl[D]\) and \(d_1\) is a conjugate over \(D\), then 
\(d_p\phi(d)\iff{}d_p\phi(d_1)\) for any definable set \(\phi\), since 
\(d_p\phi\) is defined over \(D\). Thus \(d_p\phi\) has a well defined truth 
value on the set \(\Wbar{d}\) of conjugates of \(d\), and the same is true of 
\(d_q\phi\). Since \(\Wbar{d}\) is in \(D\), this truth value is actually the 
same for \(p\) and \(q\).

This shows that we may assume \(D\) to be algebraically closed. But now 
\(d_p\) and \(d_q\) are both definition schemes for the type \(r\), so by 
uniqueness \(d_p=d_q\), and therefore \(p=q\).

The crucial point in the above argument is that \(C_a\) contains (the 
parameters for) all the defining formulas \(d_p\phi\). The definable closure 
of all these formulas is called the \Def{canonical base} of \(p\), denoted 
\(\Cb{p}\). In this terms, the previous paragraph proves the following 
proposition in the case that \(C\) is algebraically closed:

\begin{prop} \label{prp:stability}
  Let \(T\) be a stable theory with EI. There is a mapping assigning to 
  any type \(p\) over a (definably closed) set \(C\), a subset 
  \(\Cb{p}\subseteq{}C\) such that the following holds:
  \begin{enumerate}
    \item
      For any automorphism \(\tau\) (of a saturated model \(M\supseteq{}C\)),
      \begin{equation*}
        \Cb{\tau(p)}=\tau(\Cb{p})
      \end{equation*}

    \item
      For any automorphism \(\tau\) fixing \(C\) (as a set), \(\tau\) fixes 
      \(p\) if and  only if \(\tau\) fixes \(\Cb{p}\) pointwise.  
      Furthermore, \(\Cb{p}\) is the set of elements fixed by all such 
      automorphisms.

    \item
      (Existence) If \(D\supseteq{}C\) and 
      \(\Acl[\Cb{p}]\cap{}D\subseteq{}C\), then there is a type \(p_1\) over 
      \(D\) extending \(p\) with \(\Cb{p_1}=\Cb{p}\) (note that if \(C\) is 
      algebraically closed, any \(D\) satisfies the condition.)

    \item
      (Uniqueness) For any set \(A\subseteq{}C\) containing \(\Cb{p}\), \(p\) 
      is unique among extensions \(q\) of \(\Res{p}{A}\) to \(C\) with 
      \(\Cb{q}\subseteq{}A\).
  \end{enumerate}

  If \(D\) is any set such that \(D\cap\Acl[\Cb{p}]=\Cb{p}\), we denote by 
  \(\Par{p}{D}\) the restriction to \(D\) of the (unique) extension of \(p\) 
  to \(C\cup{}D\).
\end{prop}

A similar result appears in~\Cite{Bouscaren}. As for algebraically closed 
\(C\), this result implies claim~\ref{clm:cb}.

\begin{proof}
  The case when \(C\) is algebraically closed was explained above. The 
  general  statement is proved by reducing to this case.

  Let \(q\) be an extension of \(p\) to \(\Acl[C]\), and for any 
  \(b\in\Cb{q}\subseteq\Acl[C]\) let \(\Wbar{b}\) be the set of conjugates of 
  \(b\) over \(C\). We set \(\Cb{p}=B\df\{\Wbar{b}\|b\in\Cb{q}\}\). Since 
  \(C\) is definably closed, \(B{\subseteq}{}C\).
  \begin{enumerate}
    \item
      This is obvious from the definition and the algebraically closed case.

    \item
      Let \(\tau\) be an automorphism of \(M\) fixing \(p\). Then 
      \(q_1=\tau(q)\) also extends \(p\). Therefore, there is an automorphism 
      \(\sg\) of \(M\) fixing \(C\) pointwise and taking \(q_1\) to \(q\).  
      Hence \(\sg(\tau(q))=q\), so by the algebraically closed case, 
      \(\sg\circ\tau\) fixes \(\Cb{q}\) pointwise. Since \(\sg\) fixes \(C\) 
      it takes elements of \(\Acl[C]\) to their conjugates, hence so does 
      \(\tau\). Therefore \(\tau\) fixes \(\Cb{p}\), as required. The 
      furthermore part follows by the same argument.

      Conversely, if \(\tau\) fixes \(\Cb{p}\) pointwise, it takes \(\Cb{q}\) 
      to a conjugate over \(C\). Let \(\sg\) be an automorphism fixing \(C\) 
      pointwise, such that \(\sg\circ\tau\) fixes \(\Cb{q}\) pointwise. If 
      \(\tau\) preserves \(C\) as a set, it also preserves \(\Acl[C]\).  
      Therefore, by the algebraically closed case, \(\sg\circ\tau\) fixes 
      \(q\). This shows that \(\tau\) takes any extension of \(p\) to 
      \(\Acl[C]\) to another such extension. Since \(p\) is the restriction 
      to \(C\) of the intersection of all such extensions, \(\tau\) fixes 
      \(p\).

    \item
      Let \(D\) be as in the assumption. Using existence for the 
      algebraically closed case, we find \(q_1\) over \(\Acl[D]\) extending 
      \(q\), and set \(p_1=\Res{q_1}{D}\). Then any element of \(\Cb{p_1}\), 
      being a finite set of elements algebraic over \(\Cb{p}\), is itself 
      algebraic over \(\Cb{p}\).  Since it is also in \(D\), by the 
      assumption on \(D\) we have \(\Cb{p_1}\subseteq{}C\). Hence 
      \(\Cb{p_1}=\Cb{p}\).

    \item
      Assume that \(p_1\) is another extension of \(\Res{p}{A}\) to \(C\) 
      such that \(\Cb{p_1}\subseteq{}A\). Let \(q_1\) be an extension of 
      \(p_1\) to \(\Acl[C]\), and let \(r\) and \(r_1\) be the restrictions 
      of \(q\) and \(q_1\) to \(\Acl[A]\). Since \(r\) and \(r_1\) have the 
      same restriction to \(A\), there is an automorphism \(\tau\) over \(A\) 
      taking \(r\) to \(r_1\).

      However, by existence and uniqueness in the algebraically closed case, 
      \(\Cb{r}=\Cb{q}\), so \(\tau\) takes \(\Cb{r}\) to a conjugate over 
      \(C\).
      
      Therefore, as above, we may find an automorphism \(\sg\) over \(C\), 
      such that \(\sg\circ\tau\) fixes \(\Cb{r}\). Since both automorphisms 
      fix \(A\) pointwise, they fix \(\Acl[A]\) as a set, and so \(\sg\) 
      takes \(r_1\) to \(r\). Since \(\sg\) fixes \(C\), it takes \(q_1\) to 
      a type over \(\Acl[C]\) whose restriction to \(\Acl[A]\) is \(r\). By 
      the algebraically closed case, \(\sg(q_1)=q\). Since \(\sg\) fixes 
      \(C\) this means that \(p_1=p\).\qedhere
  \end{enumerate}
\end{proof}

\begin{remark}\label{rmk:cb}
  We summarise the relations between the various sets we considered. Assume 
  in claim~\ref{clm:cb} that the model \(M\) is saturated. We claim that 
  \(C_A\) is the canonical base of \(p_A=\Tp{A}{C}\). Indeed, let \(\tau\) be 
  an automorphism of \(M\) fixing \(A\) pointwise and \(C\) as a set.  Then 
  \(p_A\) is fixed by \(\tau\). By the second property above, \(\tau\) fixes 
  \(\Cb{p_A}\) pointwise. This shows that \(\Cb{p_A}\subseteq{}C_A\).  
  Conversely, let \(x\in{}C_A\), and let \(\tau\) be an automorphism fixing 
  \(C\) as a set and fixing \(p_A\). Since \(\tau\) fixes \(p_A\), \(A\) and 
  \(\tau(A)\) have the same type over \(C\). Let \(\sg\) be an automorphism 
  over \(C\) taking \(\tau(A)\) to \(A\). Then \(\sg\circ\tau\) fixes \(A\) 
  pointwise and \(C\) as a set. Hence \(\sg(\tau(x))=x\). Since \(x\in{}C\) 
  and \(\sg\) fixes \(C\), this implies that \(\tau(x)=x\). Again by the 
  second property, this shows that \(x\in\Cb{p_A}\).

  Hence, if the conditions of claim~\ref{clm:cb} hold (and \(M\) is 
  saturated), \(C_A=C_B\), and both are equal to the canonical base of the 
  type.  Considering proposition~\ref{prp:stable} again, we see that this set 
  also coincides (under the conditions of the proposition) with \(\Cc(B_i)\), 
  where \(B_i=\Dcl(B\cup{}a_i)\).
\end{remark}\index{stability|)}


\bibliographystyle{amsplain}

\bibliography{../defaut,../../../bibtex/mr}

\end{document}